\documentclass[10pt]{amsart}
\usepackage{amssymb}

\usepackage{amsfonts,amssymb}
\usepackage{amsmath,amscd}
\usepackage[dvips]{graphicx}

\usepackage[all]{xy}
\usepackage{color}

\newtheorem{Thm}{Theorem}[section]
\newtheorem{Lem}[Thm]{Lemma}
\newtheorem{Def}[Thm]{Definition}
\newtheorem{Cor}[Thm]{Corollary}
\newtheorem{Prop}[Thm]{Proposition}
\newtheorem{Ex1}[Thm]{Example}
\newtheorem{Rem1}[Thm]{Remark}

\newenvironment{Rem}{\begin{Rem1}\rm}{\end{Rem1}}

\title[Hochschild cohomology ring]{The
Batalin-Vilkovisky structure over the Hochschild cohomology ring of
a group algebra}

\author{Yuming Liu and Guodong Zhou$^*$}

\address{Yuming Liu
\newline School of Mathematical Sciences
\newline Laboratory of Mathematics and Complex Systems
\newline Beijing Normal University
\newline Beijing 100875
\newline P.R.China}
\email{ymliu@bnu.edu.cn}

\address{Guodong Zhou
 \newline Department of Mathematics
 \newline Shanghai Key laboratory of PMMP
\newline East China Normal University
\newline  Dong Chuan Road 500
\newline Shanghai 200241
 \newline P.R.China}
  \email{gdzhou@math.ecnu.edu.cn}

\date{version of \today}

\newenvironment{Proof}[1][Proof]{\begin{trivlist}
\item[\hskip \labelsep {\bfseries #1}]}{\flushright
$\Box$\end{trivlist}}

\newcommand{\lra}{\longrightarrow}

\newcommand{\ra}{\rightarrow}
\newcommand{\sdp}{\times\kern-.2em\vrule height1.1ex depth-.05ex}
\newcommand{\epi}{\lra \kern-.8em\ra}

\newcommand{\ol}{\overline}

\newcommand{\ot}{\otimes}

 \normalbaselines
\thanks{Both authors are supported by the exchange program STIC-Asie 'ESCAP' financed by the French Ministry of Foreign Affairs.  The first author is supported by NCET Program from MOE of
China and by NNSF (No.11171325, No.11331006, No.61070251).  The second author is supported by  by Shanghai Pujiang
Program (No.13PJ1402800),  by National Natural Science Foundation of China (No.11301186) and by the Doctoral Fund of Youth Scholars of Ministry of Education of China (No.20130076120001).
}

\begin{document}

\renewcommand{\thefootnote}{\alph{footnote}}
\setcounter{footnote}{-1} \footnote{ $^*$ Corresponding author.}
\renewcommand{\thefootnote}{\alph{footnote}}
\setcounter{footnote}{-1} 
\renewcommand{\thefootnote}{\alph{footnote}}
\setcounter{footnote}{-1} 

\begin{abstract} We realize explicitly the well-known additive decomposition of the Hochschild cohomology ring of a group
algebra in the elements level. As a result, we describe the cup
product, the Batalin-Vilkovisky operator and the Lie bracket in the
Hochschild cohomology ring of a group algebra.
\end{abstract}

\maketitle

\section{Introduction}

Let $k$ be a field and $G$ a finite group. Then the Hochschild
cohomology ring of the group algebra $kG$ admits an additive
decomposition:
$$HH^*(kG)\simeq \bigoplus_{x\in X}H^*(C_G(x),k)$$
where $X$ is a set of representatives of conjugacy classes of
elements of $G$ and $C_G(x)$ is the centralizer of $x\in G$. The
proof of this isomorphism can be found in \cite{Benson} or
\cite{SW1999}. The usual proof is abstract rather than giving an
explicit isomorphism. For example, one of the key steps is to use
the so-called Eckmann-Shapiro Lemma, one need to construct some
comparison maps between two projective resolutions in order to write
it down explicitly, and this is usually difficult. In \cite{SW1999},
Siegel and Witherspoon used techniques and notations from group
representation theory to interpret the above additive decomposition
explicitly. For our purpose, we need to give an explicit isomorphism
in the elements level.

A priori, the additive decomposition gives an isomorphism of graded
vector spaces. The left handed side has a graded commutative algebra
structure given by the cup product, a graded Lie algebra structure
given by the Gerstenhaber Lie bracket (\cite{Gerstenhaber1963}), and
a Batalin-Vilkovisky (BV) algebra structure given by the
$\bigtriangleup$ operator (\cite{Tradler2008}). It would be
interesting to describe these structures in terms of pieces from the
right handed side.

For graded algebra structure, it was done by Holm for abelian groups
using computations (\cite{Holm1996}), then Cibils and Solotar gave a
conceptual proof in (\cite{CS1997}). The general case was dealt with
by Siegel and Witherspoon (\cite{SW1999}), they described the cup
product formula by notations from group representation theory. Our
goal in the present paper is to represent the cup product, the Lie
bracket and the BV operator in the Hochschild cohomology ring in
terms of the additive decomposition. This is based on the explicit
construction of an isomorphism in the additive decomposition
(although there is no canonical choice for such an isomorphism).

The main obstruction in realizing an isomorphism in the additive
decomposition comes from the fact that, it is usually difficult to
construct the comparison map between two projective resolutions of
modules. There is a surprising way to simplify such construction,
namely, one can reduce it to construct a setwise self-homotopy over
one projective resolution, which is often much easier. This method
was already used in a recent paper by the second author jointly with
Le (\cite{LeZhou2011}). For convenience, we shall give a brief
introduction to this idea in Section 2.

\medskip
This article is organized as follows. In Section 2, we recall Le and
Zhou's method on constructing comparison maps. In Section 3 and 4,
we review the definitions of various structures over Hochschild
cohomology and group cohomology, using the normalized bar
resolutions. We always use the normalized bar resolutions since they
are easy to describe and can greatly simplify the computations.

In Section 5, we give a way to realize explicitly the additive
decomposition of the Hochschild cohomology of a group algebra. The
main line of our method follows from \cite{SW1999}. In Section 6, we
shall use some idea from \cite{CS1997} to give another way to
realize the additive decomposition.

We give the cup product formula in Section 7. Our formula shows that
the group cohomology $H^*(G, k)$ can be seen as a subalgebra of the
Hochschild cohomology $HH^*(kG)$ in the complex level, and that the
additive decomposition naturally gives an isomorphism of graded
$H^*(G, k)$-modules.

We deal with the $\bigtriangleup$ operator and the graded Lie
bracket in the next section. In particular, we show that the
operator $\bigtriangleup$ restricts to each summand under the
additive decomposition, and that$H^*(G, k)$ is indeed a BV
subalgebra of $HH^*(kG)$.

In the final section, we use our formulae to compute the BV structure of the Hochschild cohomology
ring for symmetric group of degree 3 over $\mathbb{F}_3$. To the best of our knowledge, this is the first concrete computation on the BV structure of
a non-commutative algebra.

\bigskip

\section{How to construct comparison
morphisms? \cite[Appendix]{LeZhou2011}}\label{Method}

\begin{Def} \rm{(cf. \cite{BianZhangZhang2009})} \label{weak} Let $A$ be an algebra over a field $k$.
Let
$$C^*: \cdots\longrightarrow
C_{n+1}\stackrel{d_{n+1}}{\longrightarrow}
C_{n}\stackrel{d_{n}}{\longrightarrow}
C_{n-1}\longrightarrow\cdots$$ be a chain complex of $A$-modules. If
there are maps (just as maps between sets) $s_n: C_n\longrightarrow
C_{n+1}$ such that $s_{n-1}d_n+d_{n+1}s_n=id_{C_n}$ for all $n$,
then the maps $\{s_n\}$ are called a setwise self-homotopy over the
complex $C^*$.
\end{Def}

\begin{Rem}
\begin{itemize}
\item[(i)]There is a setwise self-homotopy over a
complex $C^*$ of $A$-modules if and only if $C^*$ is an exact
complex, that is, $C^*$ is a zero object in the derived category
$D(ModA)$. Compare this with the usual self-homotopy, which is
equivalent to saying that $C^*$ is split exact, and hence it is a
zero object in the homotopy category $K(ModA)$.
\item[(ii)] Usually a setwise self-homotopy can be taken to be linear
maps, so it is a self-homotopy in the usual sense in the category of
complexes of $k$-vector spaces. In case that the exact complex is a
right bounded complex of $A$-$A$-bimodules, a setwise self-homotopy
can even be chosen as  homomorphisms of one-sided modules.
\end{itemize}

\end{Rem}

We will show how to use a setwise self-homotopy to construct a
comparison map. Let $M$ and $N$ be two $A$-modules, and let
$f:M\longrightarrow N$ be an $A$-module homomorphism. Suppose that
$P^*=(P_i,\partial_i)$ is a free resolution of $M$, and that
$Q^*=(Q_i,d_i)$ is a projective resolution of $N$. Suppose further
that there is a setwise self-homotopy $s=\{s_n\}$ over $Q^*$
(including $N$):

$$\xymatrix{\cdots \ar[r]^-{d_3} & Q_2\ar[r]^-{d_2}\ar[ld]^{s_2}\ar[d]^{id} & Q_1\ar[d]^{id}\ar[ld]^{s_1}\ar[r]^-{d_1} &\ar[d]^{id} Q_0\ar[ld]^{s_0}
 \ar[r]^-{d_0}& N\ar[d]^{id}\ar[ld]^{s_{-1}}\ar[r]& 0\\
\cdots \ar[r]_-{d_3} & Q_2\ar[r]_-{d_2} & Q_1\ar[r]_-{d_1} & Q_0
\ar[r]_-{d_0}& N\ar[r]& 0.}$$
For each $i\geq 0$, choose a basis $X_i$ for the free $A$-module
$P_i$ (the $i$-th term of $P^*$). We define inductively the maps
$f_i:X_i\longrightarrow Q_i$ as follows: for $x\in X_0,
f_0(x)=s_{-1}f\partial_0(x)$; for $i>1$ and for $x\in X_i,
f_i(x)=s_{i-1}f_{i-1}\partial_i(x)$. Extending $A$-linearly the maps
$f_i$ we get $A$-homomorphisms $f_i:P_i\longrightarrow Q_i$. It is
easy to verify that $\{f_i\}$ gives a chain map between the
complexes $P^*$ and $Q^*$. We illustrate the above procedure in the
following diagram:

$$\xymatrix{ P_n\ar@{.>}[dd]^{f_n}\ar[r]^{\partial_n} &
P_{n-1}\ar[d]^{f_{n-1}} & &
x\ar@{.>}[dd]^{f_n}\ar[r]^{\partial_n} & \partial_n(x)\ar[d]^{f_{n-1}}\\
&   Q_{n-1} \ar[ld]^{s_{n-1}}&& &f_{n-1}\partial_n(x) \ar[ld]^{s_{n-1}} \\
Q_n&&&
 s_{n-1}f_{n-1}\partial_n(x) & }$$

\medskip
We shall use the following standard homological fact.

\begin{Lem} \label{homotopy-equivalence} Let $A$ and $B$ be two rings and let $F: ModA\longrightarrow ModB$ be an additive contravariant (resp., covariant) functor. If $C^*$ and $D^*$ are two projective resolutions of an $A$-module $M$, then the cochain complexes $FC^*$ and $FD^*$ of $B$-modules are homotopic. In particular, if $\varphi:C^*\longrightarrow D^*$ and $\psi: D^*\longrightarrow C^*$ are two chain maps inducing identity maps $id_M: M\longrightarrow M$, then
$F\varphi: FD^*\longrightarrow FC^*$ (resp., $F\varphi:
FC^*\longrightarrow FD^*$) and $F\psi: FC^*\longrightarrow FD^*$
(resp., $F\psi: FD^*\longrightarrow FC^*$) are inverse homotopy
equivalences.
\end{Lem}

\section{Remainder on Hochschild cohomology}

In this section, we recall the definitions of various structures over Hochschild cohomology. For the cup product
and the Lie bracket in the Hochschild cohomology ring, we refer to Gerstenhaber's original paper \cite{Gerstenhaber1963}; for the
Batalin-Vilkovisky algebra structure, we refer to Tradler \cite{Tradler2008}. 

Let $k$ be a field and $A$ an associative $k$-algebra with identity
$1_A$. Denote by $\overline{A}$ the quotient space $A/(k\cdot1_A)$.
We shall write $\otimes$ for $\otimes_k$ and $A^{\otimes n}$ for the
$n$-fold tensor product $A\otimes \cdots \otimes A$. The {\it
normalized bar resolution} $(Bar_*(A),d_*)$ of $A$ is a free
resolution of $A$ as $A$-$A$-bimodules, where
$$Bar_{-1}(A)=A, \text{ and for }n\geq 0, \quad Bar_n(A)=A\otimes \overline{A}^{\otimes n}\otimes A,$$
$$d_0: Bar_0(A)=A\otimes A\longrightarrow A, \quad a_0\otimes a_1\longmapsto a_0a_1 (\text{multiplication map}), \text{ and for }n\geq 1,$$
$$d_n: Bar_n(A)\longrightarrow Bar_{n-1}(A), \quad a_0\otimes \overline{a_1}\otimes \cdots \otimes\overline{a_n}\otimes a_{n+1}\longmapsto
\sum_{i=0}^{n}(-1)^ia_0\otimes \cdots \otimes \overline{a_ia_{i+1}}\otimes
\cdots \otimes a_{n+1}.$$ The normalized bar resolution is a
natural quotient complex of the usual bar resolution. The
exactness of the normalized bar resolution is an easy consequence of
the following fact: there is a setwise self-homotopy $s_n:
Bar_n(A)\longrightarrow Bar_{n+1}(A)$ over $Bar_*(A)$ given by
$$s_n(a_0\otimes \overline{a_1}\otimes \cdots \otimes\overline{a_n}\otimes a_{n+1})= 1\otimes
a_0\otimes \overline{a_1}\otimes \cdots \otimes\overline{a_n}\otimes a_{n+1}.$$ Notice that here
each $s_n$ is just a right $A$-module homomorphism. For simplicity, in the following we will write $a_i$ for $\overline{a_i}$.

\medskip
Let $_AM_A$ be an $A$-$A$-bimodule. Remember that any
$A$-$A$-bimodule can be identified with a left module over the
enveloping algebra $A^e=A\otimes A^{op}$. We have the {\it
Hochschild cohomology complex} $(C^*(A,M),\delta_*)$:
$$C^n(A,M)=Hom_{A^e}(Bar_n(A),M)\simeq Hom_k(\overline{A}^{\otimes n},M), \quad \text{for }n\geq 0,$$
$$\delta_n: C^n(A,M)\longrightarrow C^{n+1}(A,M), \quad f\longmapsto
\delta_n(f), \quad \text{where }\delta_n(f) \text{ sends }
a_1\otimes \cdots \otimes a_{n+1} \text{ to }$$
$$a_1f(a_2\otimes \cdots \otimes
a_{n+1})+\sum_{i=1}^{n}(-1)^if(a_1\otimes \cdots \otimes
a_ia_{i+1} \otimes \cdots \otimes a_{n+1})+(-1)^{n+1}f(a_1\otimes
\cdots \otimes a_n)a_{n+1}.$$ For $n\geq 0$, the {\it degree-n
Hochschild cohomology group} of the algebra $A$ with coefficients
in $M$ is defined to be
$$HH^n(A,M)=H^n(C^*(A,M))\simeq Ext_{A^e}^n(A,M).$$

If in particular, $A=kG$ the group algebra of a finite group $G$, then  the
Hochschild cohomology complex  $(C^*(A,M),\delta_*)$ has the following form:
$$C^n(kG, M) \simeq  Hom_k(\overline{kG}^{\otimes n},M)\simeq Map(\overline{G}^n, M), \quad \text{for }n\geq 0,$$
where $\overline{G}=G-\{1\}$ and
$Map(\overline{G}^{\times n},M)$ denotes all the maps between the
sets $\overline{G}^{\times n}$ and $M$,  and the differential is given by
$$\delta_n: Map(\overline{G}^{\times n},M)\longrightarrow Map(\overline{G}^{\times n+1},M), \quad f\longmapsto
\delta_n(f),$$  where $\delta_n(f)$   sends $(g_1, \cdots,
g_{n+1})\in \overline{G}^{n+1}$   to
$$g_1f(g_2,  \cdots,
g_{n+1})+\sum_{i=1}^{n}(-1)^if(g_1,  \cdots,  g_ig_{i+1},  \cdots,  g_{n+1})+(-1)^{n+1}f(g_1,  \cdots,
g_n)g_{n+1}.$$

When $M=A$ with the obvious $A$-$A$-bimodule structure, we write $C^n(A)$ (resp.
$HH^n(A)$) for $C^n(A,A)$ (resp. $HH^n(A,A)$). Let $f\in C^n(A)$,
$g\in C^m(A)$. Then the {\it cup product} $f\cup g\in C^{n+m}(A)$ is
defined as follows:
$$f\cup g: \overline{A}^{\otimes (n+m)}\longrightarrow A,  \quad a_1\otimes \cdots \otimes a_{n+m}\longmapsto f(a_1\otimes
 \cdots \otimes a_n)g(a_{n+1}\otimes
 \cdots\otimes a_{n+m}).$$
 This cup product is associative and induces a well-defined product over $$HH^*(A)=\bigoplus_{n\geq 0}HH^n(A)=\bigoplus_{n\geq 0}Ext_{A^e}^n(A,A),$$
 which is called the {\it Hochschild cohomology ring} of $A$.
 Moreover, $HH^*(A)$ is graded commutative, that is, $\alpha\cup \beta=(-1)^{mn}\beta\cup \alpha$ for $\alpha\in HH^n(A)$ and $\beta\in
 HH^m(A)$. As usual, we call an element $\alpha\in HH^n(A)$ homogeneous of degree $n$, and its degree will be denoted by
 $|\alpha|$.

 The Lie bracket is defined as follows. Let $f\in C^n(A,M)$,
$g\in C^m(A)$. If $n,m\geq 1$, then for $1\leq i\leq n$, the
so-called {\it brace operation} $f\circ_i g\in C^{n+m-1}(A,M)$ is
defined by
$$f\circ_i g(a_1\otimes
 \cdots a_{n+m-1})=f(a_1\otimes
 \cdots \otimes a_{i-1}\otimes g(a_i\otimes \cdots \otimes a_{i+m-1})\otimes a_{i+m}\otimes \cdots \otimes a_{n+m-1});$$
if $n\geq 1$ and $m=0$, then $g\in A$ and for $1\leq i\leq n$, set
$$f\circ_i g(a_1\otimes
 \cdots a_{n-1})=f(a_1\otimes
 \cdots \otimes a_{i-1}\otimes g\otimes a_i\otimes \cdots \otimes a_{n-1});$$
 for any other case, set $f\circ_i
g$ to be zero. Define
$$f\circ g=\sum_{i=1}^n(-1)^{(m-1)(i-1)}f\circ_i g\in C^{n+m-1}(A,M)$$ and
for $f\in C^n(A)$, $g\in C^m(A)$, define $$[f, g]=f\circ
g-(-1)^{(n-1)(m-1)}g\circ f\in C^{n+m-1}(A).$$ The above $[~,~]$
induces a well-defined {\it (graded) Lie bracket} in Hochschild
cohomology
$$[~,~]: HH^n(A)\times HH^m(A)\longrightarrow HH^{n+m-1}(A)$$
such that $(HH^*(A), \cup, [~,~])$ is a {\it Gerstenhaber algebra},
that is,  for homogeneous elements $\alpha,\beta,\gamma$ in
$HH^*(A)$, the following three conditions hold:
\begin{itemize}
\item  $(HH^*(A), \cup)$ is an associative algebra and it is graded
commutative, that is, the cup product $\cup$ is an associative
multiplication and satisfies $\alpha\cup
\beta=(-1)^{|\alpha||\beta|}\beta\cup \alpha$;

\item  $(HH^*(A), [~,~])$ is a graded
Lie algebra, that is, the Lie bracket $[~,~]$ satisfies $[\alpha,\beta]=-(-1)^{(|\alpha|-1)(|\beta|-1)}[\beta,\alpha]$ and the graded Jacobi identity;

\item Possion rule: $[\alpha\cup \beta, \gamma]=[\alpha, \gamma]\cup \beta + (-1)^{|\alpha|(|\gamma|-1)}\alpha\cup[\beta,
\gamma]$.
\end{itemize}

We now assume that $A$ is a {\it symmetric $k$-algebra}, that is,
$A$ is isomorphic to its dual $D(A)=Hom_k(A,k)$ as $A^e$-modules, or
equivalently, if there exists a symmetric associative non-degenerate
bilinear form $\langle~,~\rangle: A\times A\longrightarrow k$. This
bilinear form induces a duality between the Hochschild cohomology
and the Hochschild homology. In fact, for any $n\geq 0$ there is an
isomorphism between $HH^n(A)$ and $HH_n(A)$ induced by the following
canonical isomorphisms
$$\mathrm{Hom}_k(A\otimes_{A^e}\mathrm{Bar}_n(A), k)\simeq
\mathrm{Hom}_{A^e}(\mathrm{Bar}_n(A), D(A))\simeq
\mathrm{Hom}_{A^e}(\mathrm{Bar}_n(A), A).$$ Via this duality, we
have, for $n\geq 1$, an operator $\bigtriangleup:
C^n(A)\longrightarrow C^{n-1}(A)$ which corresponds to the {\it
Connes' $B$-operator} (denoted by $B$) on the Hochschild homology
complex. More precisely, for any  $f\in C^n(A)$,
$\bigtriangleup(f)\in C^{n-1}(A)$ is given by the equation
$$\langle\bigtriangleup(f)(a_1\otimes
 \cdots \otimes a_{n-1}),a_n\rangle=\sum_{i=1}^n(-1)^{i(n-1)}\langle f(a_i\otimes
 \cdots \otimes a_{n-1}\otimes a_n\otimes a_1\otimes
 \cdots \otimes a_{i-1}),1\rangle.$$
From the well-known properties of the Connes' $B$-operator $B$ (cf.
\cite[Chapter 2]{Loday}), it is easy to see that the operator
$\bigtriangleup$ is a chain map such that the induced operation
$\bigtriangleup$ on Hochschild cohomology $HH^*(A)$ squares to zero
(in fact, $\bigtriangleup^2=0$ holds on normalized Hochschild
cochain complex level). It turns out that the Gerstenhaber algebra
$(HH^*(A), \cup, [~,~])$ together with the operator $\bigtriangleup$
is a {\it Batalin-Vilkovsky algebra} (BV-algebra), that is, in
addition to be a Gerstenhaber algebra, $(HH^*(A),\bigtriangleup)$ is
a complex and
$$[\alpha,\beta]=-(-1)^{(|\alpha|-1)|\beta|}(\bigtriangleup(\alpha\cup\beta)-\bigtriangleup(\alpha)\cup\beta-(-1)^{|\alpha|}\alpha\cup\bigtriangleup(\beta))$$
for all homogeneous elements $\alpha, \beta\in HH^*(A)$.

\begin{Rem}\label{sign-convention} The sign in the
definition of a BV-algebra depends on the choice of the definitions
of cup product and Lie bracket. If we define $\alpha{\cup}
'\beta=(-1)^{|\alpha||\beta|}\alpha\cup\beta$ and
${\bigtriangleup}'(\alpha)=(-1)^{(|\alpha|-1)}\bigtriangleup(\alpha)$,
then we get
$$[\alpha,\beta]=(-1)^{|\alpha|}({\bigtriangleup}'(\alpha{\cup}'\beta)-
{\bigtriangleup}'(\alpha){\cup}'\beta-(-1)^{|\alpha|}\alpha{\cup}'{\bigtriangleup}'(\beta)),$$
which is the equality in the usual definition of a BV-algebra (see,
for example \cite{Getzler1994,Menichi2004}). We choose the sign
convention from \cite{Tradler2008} because of our convention of the
definitions of cup product and Connes' $B$-operator in the
Hochschild (co)homology theory.
\end{Rem}

\section{Remainder on group cohomology}

Let $G$ be a finite group and $U$ a left $kG$-module. The group
cohomology of $G$ with coefficient in $U$ is defined to be $H^n(G,
U)=\mathrm{Ext}^n_{kG}(k, U)$. The complex $Bar_*(kG)\otimes_{kG}
k$  is the  standard resolution of the trivial module $k$.   In
fact, as the  setwise self-homotopy $s_n$ over  $Bar_*(kG)$ are right
module homomorphisms, $Bar_*(kG)\otimes_{kG}  k$ is exact and thus
a projective resolution of $kG\otimes_{kG}k\simeq k$. We write the
complex $C^*(G, U)=Hom_{kG}(Bar_*(kG)\otimes_{kG}  k, U)$.
Therefore,  for $n\geq 0$,
$$\begin{array}{rcl}C^n(G,U)&\simeq&  Hom_{kG}((kG\otimes \overline{kG}^{\otimes n}\otimes kG)\otimes_{kG}k, U)\simeq Hom_{kG}(kG\otimes \overline{kG}^{\otimes n}, U)\\&\simeq &Hom_{k}( \overline{kG}^{\otimes n}, U)\simeq
Map(\overline{G}^{\times n},U),\end{array}$$  and the differential
is given by
 $$\delta_0(x)(g)=gx-x  \quad (\mbox{for } x\in U \mbox{ and }
g\in \overline{G})$$ and (for $\varphi: \overline{G}^{\times
n}\longrightarrow U$ and $g_1,\cdots, g_{n+1}\in \overline{G}$)
$$\quad  \quad \quad  \delta_n(\varphi)(g_1,\cdots, g_{n+1})=g_1\varphi(g_2,\cdots,
g_{n+1}) + \quad  \quad  \quad  \quad  \quad  \quad  \quad
\quad \quad  \quad  \quad  \quad  \quad  \quad  \quad  \quad
\quad \quad \quad  \quad \quad  \quad  \quad  \quad  \quad  \quad$$
$$\quad  \quad  \quad  \quad  \quad  \quad  \quad
\quad \quad  \quad  \quad  \quad  \quad  \quad  \quad  \quad
\sum_{i=1}^n(-1)^i\varphi(g_1,\cdots, g_ig_{i+1},\cdots,
g_{n+1})+(-1)^{n+1}\varphi(g_1,\cdots, g_n).$$


Of particular interest to us are the following two
cases which relate  group cohomology to Hochschild cohomology and
 in fact which underly our two realisations of
  the additive decomposition of the Hochschild cohomology of a  group algebra.

   Note that we have an algebra isomorphism
$(kG)^e\simeq k(G\times G)$ given by $g_1\otimes g_2\longmapsto
(g_1,g_2^{-1})$, for $g_1,g_2\in G$. Thus we can also identify
each $kG$-$kG$-bimodule
 $M$ as a left $k(G\times G)$-module by $(g_1,g_2)\cdot x=g_1xg_2^{-1}$. In the sequel, we shall write the
Hochschild cohomology complex for the group algebra $kG$ in terms of
$k(G\times G)$-modules.

\textbf{Case 1}. $M=kG$, the  module $kG$ with the obvious
$kG$-$kG$-bimodule, or equivalently, the $k(G\times G)$-module
$kG$ with action: $(g_1,g_2)\cdot x=g_1xg_2^{-1}$ for $g_1,g_2\in
G$. Consider $G$ as a subgroup of $G\times G$ via the diagonal
embedding $G\to G\times G, g\mapsto (g, g)$, and it is easy to
verify that there is a $k(G\times G)$-module isomorphism
$Ind_G^{G\times G} k=k(G\times G)\otimes_{kG}k\simeq kG$,
$(g_1,g_2)\otimes 1\longmapsto g_1g_2^{-1}$. So we have
$$\begin{array}{rcl}HH^n(kG, kG)&\simeq& Ext_{k(G\times G)}^n(kG,kG)\simeq Ext_{k(G\times
G)}^n(Ind_G^{G\times G} k, kG)\\
&\simeq & Ext_{kG}^n(k, Res^{G\times G}_GkG)=Ext_{kG}^n(k,  {}_ckG)\\
&=&H^n(G, {}_ckG),\end{array}$$ where the third  isomorphism is
given by the adjoint equivalence and ${}_ckG$ is considered as a
left $kG$-module by conjugation: $g\cdot x=gxg^{-1}$ for $g, x\in
G$. This verifies a well-known fact observed by Eilenberg and Mac
Lane (\cite{EM1947}): the Hochschild cohomology $HH^n(kG,kG)$ of
$kG$ with coefficients in $kG$ is isomorphic to the ordinary group
cohomology $H^n(G,kG)$ of $G$ with coefficients in $kG$ under the
conjugation.

\textbf{Case 2.} $M=k$, the trivial $kG$-$kG$-bimodule, or
equivalently, the $k(G\times G)$-module $k$ with action:
$(g_1,g_2)\cdot 1=1$ for $g_1,g_2\in G$. Since we have
$$HH^n(kG, k)\simeq Ext_{k(G\times G)}^n(kG,k)\simeq Ext_{k(G\times G)}^n(k(G\times
G)\otimes_{kG}k,k)\simeq Ext_{kG}^n(k,k)=H^n(G, k),$$ the
Hochschild cohomology $HH^n(kG,k)$ of $kG$ with coefficients in
$k$ is isomorphic to the ordinary group cohomology $H^n(G,k)$.
Another way to see this lies in the fact that the two complexes
$C^*(kG, k)$ and $C^*(G, k)$ coincide.

We can deduce the second case from the first one. In fact, the
subspace $k(\sum_{g\in G}g)\subseteq kG$ is a sub-$(G\times
G)$-module of $kG$  (and also sub-$G$-module of ${}_ckG$),  which
is isomorphic to the trivial module. Via the isomorphisms in Case
1, $HH^*(kG, k(\sum_{g\in G}g))$ corresponds to $H^*(G,
k(\sum_{g\in G}g))$.


We can in fact define a cup product and  Lie bracket  over
$H^*(G,k)=\bigoplus_{n\geq 0}H^n(G,k)$ such that it becomes a Gerstenhaber algebra. One sees
that the cup product and the Lie bracket over $HH^*(kG)$ restrict
to $H^*(G,k)$  by \cite[Corollary 2.2]{FS2004}, so $H^*(G,k)$ is a
Gerstenhaber subalgebra of $HH^*(kG)$. In fact, as in \cite[Proof
of Theorem 1.8]{FS2004}, there is a chain map at the cohomology
complex level:
$$Hom_{kG}(Bar_n(kG)\otimes_{kG}k,k)=C^n(kG,k)\hookrightarrow C^n(kG)=Hom_{k(G\times
G)}(Bar_n(kG),kG),$$
$$(\varphi: \overline{G}^{\times n}\longrightarrow k)\longmapsto (\psi: \overline{G}^{\times n}\longrightarrow kG),  \quad \psi(g_1,
\cdots, g_{n})=\varphi(g_1, \cdots, g_{n})g_1\cdots g_{n}.$$ This
inclusion map preserves the brace operations in the following sense:

Let $\varphi_1\in C^n(kG,k)\simeq
Map(\overline{G}^{\times n},k), \varphi_2\in C^m(kG,k)$, and
let $\widehat{\varphi}_1\in C^n(kG),\widehat{\varphi}_2\in C^m(kG)$ be the
corresponding elements under the above inclusion map. Then
$\widehat{\varphi}_1\circ_i
\widehat{\varphi}_2=\widehat{\varphi_1\circ_i
\widehat{\varphi}_2}\in C^{m+n-1}(kG)$.

Recall that $kG$ is a symmetric algebra with the bilinear form
$$\langle~,~\rangle: kG\times kG\longrightarrow k,$$
$$\langle g,h\rangle=\left\{\begin{array}{cc}
1 & \mbox{ if }g=h^{-1} \\ 0 & \mbox{ otherwise }
\end{array}\right. $$ for $g,h\in G$.
So there is a well-defined BV-algebra structure on $HH^*(kG)$. We
shall see later that $H^*(G,k)$  is furthermore  a sub-BV-algebra of
$HH^*(kG)$.


\section{The first realization of the additive decomposition}
Let $k$ be a field and $G$ a finite group. Then the Hochschild
cohomology ring of the group algebra $kG$ admits an additive
decomposition:
$$HH^*(kG)\simeq \bigoplus_{x\in X}H^*(C_G(x),k),$$
where $X$ is a set of representatives of conjugacy classes of
elements of $G$ and $C_G(x)=\{g\in G\mid gx=xg\}$ is the centralizer
subgroup of $G$. In this section, we give an explicit construction
of the additive decomposition. The main technique we used here is to
construct comparison maps based on some setwise self-homotopys.

The following is a proof of the additive decomposition which
consists of a series of isomorphisms. Our first realization of the
additive decomposition will follow this series of isomorphisms.
$$\begin{array}{rcl}&&HH^*(kG, kG)
=Ext^*_{(kG)^e}(kG, kG) \simeq Ext^*_{k(G\times G)}(kG, kG)\\

&\stackrel{(1)}{\simeq} & Ext^*_{k(G\times G)}(Ind_G^{G\times G}
k,
 kG) \\

&& \mathrm{because}\  {}_{k(G\times G)}kG\simeq Ind_G^{G\times G}
k={}_{k(G\times G)}k(G\times G)\otimes_{kG}k \ \mathrm{where}\
 k(G\times G)  \ \mathrm{is\  endowed\  with\  }  \\
 &&\mathrm{the \ right\ }kG\mathrm{-module\ structure\ }\ \mathrm{via\  the\  diagonal\  map\ } G\to
G\times G, g\mapsto
(g, g)\\

   &\stackrel{(2)}{\simeq} &
Ext^*_{kG}(k, Res_G^{G\times G}kG)=Ext^*_{kG}(k, Hom_{k(G\times
G)}(k(G\times
G), kG))=Ext_{kG}^*(k, {}_ckG)=\\
&& H^*(G, {}_ckG)\mathrm{\ by \ the  \ adjoint \ pair \
({}_{k(G\times G)}k(G\times G)\otimes_{kG}-,Hom_{k(G\times
G)}({}_{k(G\times G)}k(G\times G)_{kG},-)) } ¡¡
  \\

   &\stackrel{(3)}{\simeq} &\oplus_{x\in X} Ext^*_{kG}(k, {}_ckC_x)\mathrm{\ because}\ {}_ckG=\oplus {}_{x\in X}\ {}_ckC_x
   \mathrm{\ where}{}_ckC_x\mathrm{\ is\ the\ }\\
 &&kG\mathrm{-module\ generated \ by \ the \
elements\ in \ the \ conjugacy\ class\ }C_x=\{gxg^{-1}|g\in G\} ¡¡\\

  &\stackrel{(4)}{\simeq} &\oplus_{x\in X} Ext^*_{kG}(k,
Coind_{C_G(x)}^G k)\\

&&\mathrm{because}\ \mathrm{as\ left}\ kG\mathrm{-modules},\
{}_ckC_x\simeq Coind_{C_G(x)}^G k=Hom_{kC_G(x)}(kG, k)\\

 & \stackrel{(5)}{\simeq} &\oplus_{x\in X} Ext^*_{kC_G(x)}(
Res_{C_G(x)}^G k, k)  \\
&& \mathrm{by \ the \  adjoint\ pair\ ({}_{kC_G(x)}kG\ot_{kG}-,\ Hom_{kC_G(x)}({}_{kC_G(x)}kG_{kG}, -))}\\

& \stackrel{(6)}{\simeq} &\oplus_{x\in X} Ext^*_{kC_G(x)}(k,
k))=\oplus_{x\in X}H^*(C_G(x), k)

\end{array}$$

We shall express explicitly    these isomorphisms step by step
using the bar resolution.

\bigskip

{\it The first step.} By definition, the Hochschild cohomology
groups $HH^*(kG, kG)$ can be computed using the bar resolution
$Bar_*(kG)$.  On the other hand, $Bar_*(kG)\otimes_{kG}k$ is a free
resolution of $k$ as left $kG$-module and therefore,
 $k(G\times
G)\otimes_{kG}Bar_*(kG)\otimes_{kG}k$ is also a free resolution of
the   $k(G\times G)$-module $k(G\times G)\ot_{kG}k\simeq kG$. Notice
that the terms in $Bar_*(kG)$ are still viewed as the usual
$kG$-$kG$-bimodules when we do the above tensor products.

Let us write explicitly the resolution $k(G\times
G)\otimes_{kG}Bar_*(kG)\otimes_{kG}k$. Under the identification
$$\begin{array}{rcl}k(G\times G)\otimes_{kG}Bar_n(kG)\otimes_{kG}k
&\simeq & k(G\times G)\ot_{kG}(kG\ot  \ol{kG}^{\ot n}\ot
kG)\ot_{kG} k\\
&\simeq & k(G\times G)\ot \ol{kG}^{\ot n}\\
&\simeq& kG\otimes kG\otimes \overline{kG}^{\otimes n},\end{array}
 $$
 and the
differential is as follows (we only write down the maps on base
elements here and later):
$$kG\otimes kG \longrightarrow kG, \quad x\otimes y\longmapsto
xy^{-1};$$
$$kG\otimes kG\otimes \overline{kG} \longrightarrow kG\otimes kG, \quad x\otimes
y\otimes g_1\longmapsto xg_1\otimes yg_1-x\otimes y;$$
$$\cdots \cdots \cdots$$
$$kG\otimes kG\otimes \overline{kG}^{\otimes
n}\longrightarrow kG\otimes kG\otimes \overline{kG}^{\otimes {n-1}},
 \quad x\otimes y\otimes g_1\otimes \cdots \otimes g_{n}\longmapsto$$
$$
xg_1\otimes yg_1\otimes g_2\otimes \cdots \otimes
g_{n}+\sum_{i=1}^{n-1}(-1)^ix\otimes y\otimes g_1\otimes \cdots
\otimes g_ig_{i+1}\otimes \cdots \otimes g_{n}+(-1)^nx\otimes
y\otimes g_1\otimes \cdots \otimes g_{n-1}.$$

We also have
$$Hom_{k(G\times
G)}(k(G\times G)\otimes_{kG}Bar_n(kG)\otimes_{kG}k,kG)\simeq
Hom_k(\overline{kG}^{\otimes n},kG)\simeq Map(\overline{G}^{\times
n},kG).$$ Using this identification, $H^*(Hom_{k(G\times
G)}(k(G\times G)\otimes_{kG}Bar_*(kG)\otimes_{kG}k,kG))$ is given
by the following cochain complex:
$$0\longrightarrow kG\stackrel{{\delta_0}}{\longrightarrow} Map(\overline{G},kG)\stackrel{{\delta_1}}{\longrightarrow} \cdots \longrightarrow Map(\overline{G}^{\times n},kG)
\stackrel{{\delta_n}}{\longrightarrow} \cdots,$$ where the differential
is given by $${\delta_0}(x)(g)=gxg^{-1}-x  \quad  \mbox{for } x\in kG
\mbox{ and } g\in \overline{G},$$ and  for $\varphi:
\overline{G}^{\times n}\longrightarrow kG$ and $g_1,\cdots,
g_{n+1}\in \overline{G}$,
$$\quad  \quad \quad  {\delta_n}(\varphi)(g_1,\cdots, g_{n+1})=g_1\varphi(g_2,\cdots,
g_{n+1})g_1^{-1}+ \quad  \quad  \quad  \quad  \quad  \quad  \quad
\quad \quad  \quad  \quad  \quad  \quad  \quad  \quad  \quad  \quad
\quad \quad  \quad \quad  \quad  \quad  \quad  \quad  \quad  \quad
\quad \quad \quad  \quad $$ $$\quad  \quad  \quad  \quad  \quad
\quad \quad  \quad \quad  \quad \quad  \quad  \quad \quad
\sum_{i=1}^n(-1)^i\varphi(g_1,\cdots, g_ig_{i+1},\cdots,
g_{n+1})+(-1)^{n+1}\varphi(g_1,\cdots, g_n).$$

We will show that the two complexes $k(G\times
 G)\ot_{kG}Bar_*(kG)\ot_{kG}k$ and $Bar_*(kG)$ are isomorphic and therefore there is an isomorphism
$$(1)  \quad  H^*(Hom_{k(G\times G)}(Bar_*(kG),kG))\simeq
H^*(Hom_{k(G\times G)}(k(G\times
G)\otimes_{kG}Bar_*(kG)\otimes_{kG}k,kG)).$$ To do this, we need to construct the comparison maps between the two free resolutions
$Bar_*(kG)$ and $k(G\times G)\otimes_{kG}Bar_*(kG)\otimes_{kG}k$ of
the above $k(G\times G)$-module $kG$. As explained in
Section~\ref{Method}, this is reduced to construct setwise
self-homotopys over these resolutions. Our principle here is to
choose those setwise self-homotopys so that the computations and
results are as simple as possible.

We choose a setwise self-homotopy over $Bar_*(kG)$ as follows:
$$u_{-1}: kG\to kG\otimes kG, \quad g\mapsto g\ot 1,$$ and for $n\geq
0$,
$$u_n: kG\otimes \overline{kG}^{\otimes n}\otimes kG\longrightarrow
kG\otimes \overline{kG}^{\otimes n}\otimes kG,$$
$$g_0\otimes g_1\otimes \cdots \otimes g_{n+1}\longmapsto
(-1)^{n+1}g_0\otimes g_1\otimes \cdots \otimes g_{n+1}\otimes 1.$$
Using $\{u_n\}$ we can construct a comparison map $$\alpha_*:
k(G\times G)\otimes_{kG}Bar_*(kG)\otimes_{kG}k=kG\ot kG\ot
k\ol{G}^{\ot *}\longrightarrow Bar_*(kG)=kG\ot k {G}^{\ot *}\ot
kG$$ as follows (as before we only write down the maps on base
elements):
$$\alpha_{-1}: kG\longrightarrow kG, \quad x\longmapsto x,$$
$$\alpha_{0}: kG\otimes kG\longrightarrow kG\otimes kG, \quad x\otimes y\longmapsto x\otimes y^{-1},$$
$$\alpha_{1}: kG\otimes kG\otimes \overline{kG}\longrightarrow kG\otimes \overline{kG}\otimes
kG, \quad x\otimes y\otimes g_1\longmapsto -xg_1\otimes
g_1^{-1}\otimes y^{-1},$$
$$\cdots \cdots \cdots$$
$$\alpha_{n}: kG\otimes kG\otimes \overline{kG}^{\otimes n}\longrightarrow kG\otimes \overline{kG}^{\otimes n}\otimes
kG, \quad x\otimes y\otimes g_1\otimes \cdots \otimes g_n
 \longmapsto$$
$$ (-1)^{\frac{n(n+1)}{2}}xg_1\cdots g_n\otimes g_n^{-1}\otimes \cdots
\otimes g_1^{-1}\otimes y^{-1}.$$

Similarly, we choose a setwise self-homotopy over $k(G\times
G)\otimes_{kG}Bar_*(kG)\otimes_{kG}k$ as follows:
$$ v_{-1}: kG\to kG\ot kG, \quad g\mapsto g\ot 1,$$ and for $n\geq 0$,
\begin{align*}
v_n:{} & kG\otimes kG\otimes \overline{kG}^{\otimes
n}\longrightarrow
        kG\otimes kG\otimes \overline{kG}^{\otimes n},\\
     & x\otimes y\otimes g_1\otimes \cdots \otimes g_n
        \longmapsto xy^{-1}\otimes 1\otimes y\otimes g_1\otimes \cdots
        \otimes g_n.
\end{align*}
Using $\{v_n\}$ we can construct a comparison map
$$\beta_*: Bar_*(kG)=kG\ot k {G}^{\ot *}\ot
kG\longrightarrow k(G\times
G)\otimes_{kG}Bar_*(kG)\otimes_{kG}k=kG\ot kG\ot k {G}^{\ot *}$$
as follows:
$$\beta_{-1}: kG\longrightarrow kG, \quad x\longmapsto x,$$
$$\beta_{0}: kG\otimes kG\longrightarrow kG\otimes kG, \quad x\otimes y\longmapsto x\otimes y^{-1},$$
$$\beta_{1}: kG\otimes \overline{kG}\otimes kG\longrightarrow kG\otimes
kG\otimes \overline{kG}, \quad x\otimes g_1\otimes y\longmapsto
-xg_1\otimes y^{-1}\otimes g_1^{-1},$$
$$\cdots \cdots \cdots$$
$$\beta_{n}: kG\otimes \overline{kG}^{\otimes n}\otimes kG\longrightarrow kG\otimes
kG\otimes \overline{kG}^{\otimes n}, \quad x\otimes g_1\otimes
\cdots \otimes g_n
 \otimes y\longmapsto$$
 $$ (-1)^{\frac{n(n+1)}{2}}xg_1\cdots g_n\otimes y^{-1}\otimes g_n^{-1}\otimes \cdots
\otimes g_1^{-1}.$$

It is easy to check that the chain maps $\{\alpha_n\}$ and
$\{\beta_n\}$ are inverse to each other, and therefore we get an
isomorphism $$Hom_{k(G\times
G)}(Bar_*(kG),kG)\longrightarrow Hom_{k(G\times
G)}(k(G\times G)\otimes_{kG}Bar_*(kG)\otimes_{kG}k,kG),$$
$$(\varphi: \overline{G}^{\times n}\longrightarrow kG)\longmapsto (\varphi_1:
\overline{G}^{\times n}\longrightarrow kG), \quad
\varphi_1(g_1,\cdots, g_{n})= (-1)^{\frac{n(n+1)}{2}}g_1\cdots
g_n\varphi(g_{n}^{-1},\cdots, g_1^{-1}).$$ Its inverse is given by
$$Hom_{k(G\times G)}(k(G\times
G)\otimes_{kG}Bar_*(kG)\otimes_{kG}k,kG)\longrightarrow
Hom_{k(G\times G)}(Bar_*(kG),kG),$$
$$(\varphi_1:
\overline{G}^{\times n}\longrightarrow kG)\longmapsto (\varphi:
\overline{G}^{\times n}\longrightarrow kG), \quad
\varphi(g_1,\cdots, g_{n})= (-1)^{\frac{n(n+1)}{2}}g_1\cdots
g_n\varphi_1(g_{n}^{-1},\cdots, g_1^{-1}).$$ Passing to the
cohomology, we realize an isomorphism in $(1)$ and its inverse.

\medskip

{\it The second step.} Since $({}_{k(G\times G)}k(G\times
G)\otimes_{kG}-,Hom_{k(G\times G)}({}_{k(G\times G)}k(G\times
G)_{kG},-))$ is an adjoint pair, we have an isomorphism (here
$k(G\times G)$ is viewed as a right $kG$-module by diagonal
action)
$$Hom_{k(G\times G)}(k(G\times
G)\otimes_{kG}Bar_*(kG)\otimes_{kG}k,kG)\simeq
Hom_{kG}(Bar_*(kG)\otimes_{kG}k, {}_ckG).$$   Passing to the
cohomology, we get an isomorphism $$(2)  \quad H^*(Hom_{k(G\times
G)}(k(G\times G)\otimes_{kG}Bar_*(kG)\otimes_{kG}k,kG))\simeq
H^*(Hom_{kG}(Bar_*(kG)\otimes_{kG}k, {}_ckG)).$$ Remind that the
right hand side is just the ordinary group cohomology $H^*(G,kG)$
of $G$ with coefficients in ${}_ckG$. We also have
$$Hom_{kG}(Bar_n(kG)\otimes_{kG}k,kG)\simeq
Hom_{kG}(kG\otimes\overline{kG}^{\otimes n},kG)\simeq
Hom_k(\overline{kG}^{\otimes n},kG)\simeq Map(\overline{G}^{\times
n},kG).$$ Using this identification,
$H^*(G,kG)=H^*(Hom_{kG}(Bar_*(kG)\otimes_{kG}k,kG))$ is given by
the following cochain complex:
$$0\longrightarrow kG\stackrel{\delta_0}{\longrightarrow}
Map(\overline{G},kG)\stackrel{\delta_1}{\longrightarrow} \cdots \longrightarrow Map(\overline{G}^{\times n},kG)
\stackrel{\delta_n}{\longrightarrow} \cdots,$$ where the differential is
given by $$\delta_0(x)(g)=gxg^{-1}-x  \quad (\mbox{for } x\in kG \mbox{
and } g\in \overline{G}),$$ and (for $\varphi: \overline{G}^{\times
n}\longrightarrow kG$ and $g_1,\cdots, g_{n+1}\in \overline{G}$)
$$\quad  \quad \quad  \delta_n(\varphi)(g_1,\cdots, g_{n+1})=g_1\varphi(g_2,\cdots,
g_{n+1})g_1^{-1}+ \quad  \quad  \quad  \quad  \quad  \quad  \quad
\quad \quad  \quad  \quad  \quad  \quad  \quad  \quad  \quad
\quad \quad \quad  \quad \quad  \quad  \quad  \quad  \quad  \quad
\quad \quad \quad \quad  \quad $$ $$\quad  \quad  \quad  \quad
\quad \quad \quad  \quad \quad  \quad \quad  \quad  \quad \quad
\sum_{i=1}^n(-1)^i\varphi(g_1,\cdots, g_ig_{i+1},\cdots,
g_{n+1})+(-1)^{n+1}\varphi(g_1,\cdots, g_n).$$ So formally the
left hand side and the right hand side in $(2)$ are identical,
though they have different meaning. It is also easy to check that
under the above identifications, the adjoint isomorphisms are
identity maps:
$$Hom_{k(G\times
G)}(k(G\times
G)\otimes_{kG}Bar_*(kG)\otimes_{kG}k,kG)\longrightarrow
Hom_{kG}(Bar_*(kG)\otimes_{kG}k, {}_ckG),$$
$$(\varphi_1: \overline{G}^{\times n}\longrightarrow kG)\longmapsto (\varphi_2:
\overline{G}^{\times n}\longrightarrow kG), \quad
\varphi_2(g_1,\cdots, g_{n})= \varphi_1(g_1,\cdots, g_{n}).$$ Its
inverse is given by
$$Hom_{kG}(Bar_*(kG)\otimes_{kG}k, {}_ckG)\longrightarrow
Hom_{k(G\times G)}(k(G\times
G)\otimes_{kG}Bar_*(kG)\otimes_{kG}k,kG),$$
$$(\varphi_2: \overline{G}^{\times n}\longrightarrow kG)\longmapsto (\varphi_1:
\overline{G}^{\times n}\longrightarrow kG), \quad
\varphi_1(g_1,\cdots, g_{n})= \varphi_2(g_1,\cdots, g_{n}).$$
Passing to the cohomology, we realize an isomorphism in $(2)$ and
its inverse.

\medskip {\it The third step.} We choose a complete set $X$ of
representatives of the conjugacy classes in the finite group $G$.
Take $x\in X$. Then $C_x=\{gxg^{-1}|g\in G\}$ is the conjugacy class
corresponding to $x$ and $C_G(x)=\{g\in G|gxg^{-1}=x\}$ is the
centralizer subgroup. Clearly the $k$-space $kC_x$ generated by the
elements in $C_x$ is a left $kG$-module under the conjugation
action. We choose a right coset decomposition of $C_G(x)$ in $G$:
$G=C_G(x)\gamma_{1,x}\cup C_G(x)\gamma_{2,x}\cup\cdots \cup
C_G(x)\gamma_{n_x,x}$ (equivalently, $G=\gamma_{1,x}^{-1}C_G(x)\cup
\gamma_{2,x}^{-1}C_G(x)\cup\cdots \cup \gamma_{n_x,x}^{-1}C_G(x)$ is
a left coset decomposition of $C_G(x)$ in $G$), and such that
$C_x=\{x=\gamma_{1,x}^{-1}x\gamma_{1,x},\gamma_{2,x}^{-1}x\gamma_{2,x},\cdots,
\gamma_{n_x,x}^{-1}x\gamma_{n_x,x}\}$. (We will always take
$\gamma_{1,x}=1$, and we write $x_i$ for
$\gamma_{i,x}^{-1}x\gamma_{i,x}$.) Then we have the following
$kG$-module isomorphisms:
$${}_ckC_x\simeq Ind_{C_G(x)}^Gk={}_{kG}kG\otimes_{kC_G(x)}k, \quad x_i\longmapsto \gamma_{i,x}^{-1}\otimes 1,$$
$${}_ckC_x\simeq Coind_{C_G(x)}^Gk=Hom_{kC_G(x)}({}_{kC_G(x)}kG_{kG}, {}_{kC_G(x)}k), \quad x_i\longmapsto \gamma_i: kG\longrightarrow k, \gamma_i(\gamma_{j,x})=\delta_{ij},$$
where in the first isomorphism, the left $kG$-module structure on
$kG$ is the usual left multiplication and the right
$kC_G(x)$-module structure on $kG$ is given by restriction, and
$k$ is the trivial $kC_G(x)$-module, and  the same as in the
second isomorphism.

In the second step, we have arrived at the ordinary group
cohomology $H^*(G,kG)$ of $G$ with coefficients in ${}_ckG$. This
${}_ckG$ has a $kG$-module decomposition:
$${}_ckG=\bigoplus_{x\in X}{}_ckC_x.$$
Denote by $\pi_x:kG\longrightarrow kC_x$ and
$i_x:kC_x\longrightarrow kG$ the canonical projection and the
canonical injection, respectively. Then we have the following
isomorphism
$$Hom_{kG}(Bar_*(kG)\otimes_{kG}k, {}_ckG)\longrightarrow \bigoplus_{x\in X}Hom_{kG}(Bar_*(kG)\otimes_{kG}k, {}_c kC_x),$$
$$(\varphi_2: \overline{G}^{\times n}\longrightarrow kG)\longmapsto \varphi_3=\{\varphi_{3,x}|x\in X\}, \mbox{ where }\varphi_{3,x}=\pi_x\varphi_2:
\overline{G}^{\times n}\longrightarrow kC_x.$$ Its inverse is given
by
$$\bigoplus_{x\in X}Hom_{kG}(Bar_*(kG)\otimes_{kG}k,kC_x)\longrightarrow
Hom_{kG}(Bar_*(kG)\otimes_{kG}k,kG),$$
$$\varphi_3=\{\varphi_{3,x}:\overline{G}^{\times n}\longrightarrow kC_x|x\in X\}\longmapsto
(\varphi_2=\sum_{x\in X}i_x\varphi_{3,x}: \overline{G}^{\times
n}\longrightarrow kG).$$
 Passing to the cohomology, we
realize an isomorphism: $$(3)  \quad H^*(G, {}_ckG)\simeq
\bigoplus_{x\in X}H^*(G, {}_ckC_x).$$

\medskip {\it The fourth step.} We have stated in the third step the following $kG$-module isomorphism
$${}_ckC_x\simeq Hom_{kC_G(x)}(kG,k), \quad x_i\longmapsto \gamma_i: kG\longrightarrow k, \gamma_i(\gamma_{j,x})=\delta_{ij}.$$
Therefore we have the following isomorphism
$$Hom_{kG}(Bar_*(kG)\otimes_{kG}k, {}_ckC_x)\longrightarrow Hom_{kG}(Bar_*(kG)\otimes_{kG}k,Hom_{kC_G(x)}(kG,k)),$$
$$(\varphi_{3,x}:
\overline{G}^{\times n}\longrightarrow kC_x)\longmapsto
(\varphi_{4,x}: \overline{G}^{\times n}\longrightarrow
Hom_{kC_G(x)}(kG,k)),$$ where if we write
$\varphi_{3,x}(g_1,g_2,\cdots, g_n)=\sum_{i=1}^{n_x}a_{i,x}x_i$,
then $\varphi_{4,x}(g_1,g_2,\cdots, g_n)$ maps $\gamma_{i,x}$ to
$a_{i,x}$ for any $i$. The inverse isomorphism is given by
$$Hom_{kG}(Bar_*(kG)\otimes_{kG}k,Hom_{kC_G(x)}(kG,k))\longrightarrow
Hom_{kG}(Bar_*(kG)\otimes_{kG}k, {}_ckC_x),$$
$$(\varphi_{4,x}:
\overline{G}^{\times n}\longrightarrow
Hom_{kC_G(x)}(kG,k))\longmapsto (\varphi_{3,x}: \overline{G}^{\times
n}\longrightarrow kC_x),$$ where if $\varphi_{4,x}(g_1,g_2,\cdots,
g_n)$ maps $\gamma_{i,x}$ to $a_{i,x}$ for any $i$, then
$\varphi_{3,x}(g_1,g_2,\cdots, g_n)=\sum_{i=1}^{n_x}a_{i,x}x_i$.
 Passing to the cohomology, we
realize an isomorphism: $$(4)  \quad H^*(G,kC_x)\simeq
H^*(Hom_{kG}(Bar_*(kG)\otimes_{kG}k,Hom_{kC_G(x)}(kG,k))).$$

\medskip {\it The fifth step.} Since $(kG\otimes_{kG}-,Hom_{kC_G(x)}(kG,-))$ is an adjoint pair (restriction and coinduction),
we have the following isomorphism
$$Hom_{kG}(Bar_*(kG)\otimes_{kG}k,Hom_{kC_G(x)}(kG,k))\longrightarrow Hom_{kC_G(x)}(Bar_*(kG)\otimes_{kG}k,k).$$
 Passing to the cohomology, we
get an isomorphism $$(5)  \quad
H^*(Hom_{kG}(Bar_*(kG)\otimes_{kG}k,Hom_{kC_G(x)}(kG,k)))\simeq
H^*(Hom_{kC_G(x)}(Bar_*(kG)\otimes_{kG}k,k)),$$ where the right
hand side is isomorphic to the ordinary group cohomology
$H^*(C_G(x),k)$ of $C_G(x)$ with coefficients in the trivial
module $k$. Since there are $kC_G(x)$-module isomorphisms
$$Bar_*(kG)\otimes_{kG}k\simeq \bigoplus_{i=1}^{n_x}kC_G(x)\gamma_{i,x}\otimes\overline{kG}^{\otimes *},$$
we have $$Hom_{kC_G(x)}(Bar_*(kG)\otimes_{kG}k,k)\simeq
Hom_{k}(\bigoplus_{i=1}^{n_x}k\gamma_{i,x}\otimes\overline{kG}^{\otimes
n},k)\simeq Map(S_x\times\overline{G}^{\times n},k),$$ where
$S_x=\{\gamma_{1,x},\cdots, \gamma_{n_x,x}\}$ (cf. The third
step). Using this identification, the adjoint isomorphism is given
by
$$Hom_{kG}(Bar_*(kG)\otimes_{kG}k,Hom_{kC_G(x)}(kG,k))\longrightarrow Hom_{kC_G(x)}(Bar_*(kG)\otimes_{kG}k,k),$$
$$(\varphi_{4,x}:
\overline{G}^{\times n}\longrightarrow
Hom_{kC_G(x)}(kG,k))\longmapsto (\varphi_{5,x}:
S_x\times\overline{G}^{\times n}\longrightarrow k),$$ where if
$\varphi_{4,x}(g_1,g_2,\cdots, g_n)$ maps $\gamma_{i,x}$ to
$a_{i,x}$ for any $i$, then
$\varphi_{5,x}(\gamma_{i,x},g_1,g_2,\cdots, g_n)=a_{i,x}$ for any
$i$. The inverse isomorphism is given by
$$Hom_{kC_G(x)}(Bar_*(kG)\otimes_{kG}k,k)\longrightarrow Hom_{kG}(Bar_*(kG)\otimes_{kG}k,Hom_{kC_G(x)}(kG,k)),$$
$$(\varphi_{5,x}:
S_x\times\overline{G}^{\times n}\longrightarrow k)\longmapsto
(\varphi_{4,x}: \overline{G}^{\times n}\longrightarrow
Hom_{kC_G(x)}(kG,k)),$$ where if
$\varphi_{5,x}(\gamma_{i,x},g_1,g_2,\cdots, g_n)=a_{i,x}$ for any
$i$, then $\varphi_{4,x}(g_1,g_2,\cdots, g_n)$ maps $\gamma_{i,x}$
to $a_{i,x}$ for any $i$.  Passing to the cohomology, we realize an
isomorphism in $(5)$ and its inverse.

\medskip {\it The sixth step.} In the fifth step, we have arrived at the ordinary group cohomology $H^*(C_G(x),k)$
of $C_G(x)$ with coefficients in the trivial module $k$, where
$H^*(C_G(x),k)$ is computed by the cochain complex
$Hom_{kC_G(x)}(Bar_*(kG)\otimes_{kG}k,k)$. By the identification in
fifth step, this is given by the following cochain complex:
$$0\longrightarrow k^{\times n_x}\stackrel{\delta_0}{\longrightarrow} Map(S_x\times\overline{G},k)
\stackrel{\delta_1}{\longrightarrow} \cdots \longrightarrow
Map(S_x\times\overline{G}^{\times n},k)
\stackrel{\delta_n}{\longrightarrow} \cdots,$$ where the differential is
given by $\delta_0(\{a_{i,x}\})((\gamma_{j,x},g_1))=a_{s_j,x}-a_{j,x},$
such that $a_{s_j,x}$ is determined as follows: for $\{a_{i,x}\}\in
k^{\times n_x}, \gamma_{j,x}\in S_x, g_1\in \overline{G}$, we have
$$\gamma_{j,x}g_1=h_{j,1}\gamma_{s_j,x} \mbox{ for
some } h_{j,1}\in C_G(x) \mbox{ and for some } 1\leq s_j\leq n_x,$$
and (for $\varphi: S_x\times\overline{G}^{\times
n}\longrightarrow k$, $\gamma_{j,x}\in S_x, g_1,\cdots, g_{n+1}\in
\overline{G}$ such that $\gamma_{j,x}g_1=h_{j,1}\gamma_{s_j,x}$)
$$\quad  \quad \quad  \delta_n(\varphi)(\gamma_{j,x},g_1,\cdots, g_{n+1})=\varphi(\gamma_{s_j,x},g_2,\cdots,
g_{n+1})+ \quad  \quad  \quad  \quad  \quad  \quad  \quad \quad
\quad \quad  \quad  \quad  \quad  \quad  \quad  \quad  \quad \quad
\quad \quad \quad  \quad  \quad  \quad  \quad  \quad  \quad \quad
\quad \quad  \quad $$ $$\quad  \quad  \quad  \quad
\quad \quad \quad  \quad \quad  \quad  \quad \quad
\sum_{i=1}^n(-1)^i\varphi(\gamma_{j,x},g_1,\cdots,
g_ig_{i+1},\cdots,
g_{n+1})+(-1)^{n+1}\varphi(\gamma_{j,x},g_1,\cdots, g_{n}).$$
(Remark that for a fixed $g_1\in \overline{G}$, $\{s_1,s_2,\cdots,
s_{n_x}\}$ is a permutation of $\{1,2,\cdots, n_x\}$.)

The above computation for $H^*(C_G(x),k)$ uses the projective
resolution $Bar_*(kG)\otimes_{kG}k$ of the trivial $kC_G(x)$-module
$k$, which is identified as the following complex (It is in fact a
projective resolution of the trivial $kG$-module $k$, but we view it
as a complex of $kC_G(x)$-modules by restriction)
$$\cdots \longrightarrow kG\otimes \overline{kG}^{\otimes
n}\stackrel{d_n}{\longrightarrow} \cdots \longrightarrow
kG\otimes\overline{kG}\stackrel{d_1}{\longrightarrow}
kG\stackrel{d_0}{\longrightarrow} k\longrightarrow 0,$$where the
differential is given by
$$d_0(g_0)=1  \quad (\mbox{for } g_0\in G)$$
 and (for $g_0\in G, g_1,\cdots, g_{n}\in \overline{G}$)
$$\quad  \quad \quad  d_n(g_0, g_1,\cdots, g_{n})=g_0 g_1\otimes g_2\otimes\cdots
\otimes g_{n}+ \quad  \quad  \quad  \quad  \quad  \quad \quad \quad
\quad  \quad  \quad  \quad  \quad  \quad  \quad  \quad \quad \quad
\quad  \quad \quad  \quad  \quad  \quad  \quad  \quad \quad \quad
\quad \quad  \quad $$ $$\quad  \quad  \quad  \quad \quad \quad \quad
\quad \quad  \quad \quad  \quad  \quad \quad
\sum_{i=1}^{n-1}(-1)^ig_0\otimes \cdots \otimes g_ig_{i+1}\otimes
\cdots \otimes g_{n}+(-1)^{n}g_0\otimes g_1\otimes \cdots \otimes
g_{n-1}.$$We now use another projective resolution
$Bar_*(kC_G(x))\otimes_{kC_G(x)}k$ of the trivial $kC_G(x)$-module
$k$, which is identified as the following complex
$$\cdots
\longrightarrow kC_G(x)\otimes \overline{kC_G(x)}^{\otimes
n}\stackrel{d_n}{\longrightarrow} \cdots \longrightarrow
kC_G(x)\otimes\overline{kC_G(x)}\stackrel{d_1}{\longrightarrow}
kC_G(x)\stackrel{d_0}{\longrightarrow} k\longrightarrow 0,$$ where
the differential is given by
$$d_0(h_0)=1  \quad (\mbox{for } h_0\in {C_G(x)})$$
 and (for $h_0\in C_G(x), h_1,\cdots, h_{n}\in \overline{C_G(x)}$)
$$\quad  \quad \quad  d_n(h_0, h_1,\cdots, h_{n})=h_0 h_1\otimes h_2\otimes\cdots
\otimes h_{n}+ \quad  \quad  \quad  \quad  \quad  \quad \quad \quad
\quad  \quad  \quad  \quad  \quad  \quad  \quad  \quad \quad \quad
\quad  \quad \quad  \quad  \quad  \quad  \quad  \quad \quad \quad
\quad \quad  \quad $$ $$\quad  \quad  \quad  \quad \quad \quad \quad
\quad \quad  \quad \quad  \quad  \quad \quad
\sum_{i=1}^{n-1}(-1)^ih_0\otimes \cdots \otimes h_ih_{i+1}\otimes
\cdots \otimes h_{n}+(-1)^{n}h_0\otimes h_1\otimes \cdots \otimes
h_{n-1}.$$We have
$$Hom_{kC_G(x)}(Bar_*(kC_G(x))\otimes_{kC_G(x)}k,k) \simeq
Map(\overline{C_G(x)}^{\times n},k),$$ so $H^*(C_G(x),k)$ can also
be computed by the following cochain complex
$$0\longrightarrow k\stackrel{\delta_0}{\longrightarrow} Map(\overline{C_G(x)},k)
\stackrel{\delta_1}{\longrightarrow} \cdots \longrightarrow
Map(\overline{C_G(x)}^{\times n},k) \stackrel{\delta_n}{\longrightarrow}
\cdots,$$ where the differential is given by $$\delta_0(a)(h_1)=0
  \quad  \quad  (\mbox{for } a\in k, h_1\in \overline{C_G(x)}) $$ and (for $\varphi:
\overline{C_G(x)}^{\times n}\longrightarrow k$, $h_1,\cdots,
h_{n+1}\in \overline{C_G(x)}$)
$$\quad  \quad \quad  \delta_n(\varphi)(h_1,\cdots, h_{n+1})=\varphi(h_2,\cdots,
h_{n+1})+ \quad  \quad  \quad  \quad  \quad  \quad  \quad \quad
\quad \quad  \quad  \quad  \quad  \quad  \quad  \quad  \quad \quad
\quad \quad \quad  \quad  \quad  \quad  \quad  \quad  \quad \quad
\quad \quad  \quad $$ $$\quad  \quad  \quad  \quad  \quad \quad
\quad \quad \quad  \quad \quad  \quad  \quad \quad
\sum_{i=1}^n(-1)^i\varphi(h_1,\cdots, h_ih_{i+1},\cdots,
h_{n+1})+(-1)^{n+1}\varphi(h_1,\cdots, h_{n}).$$ Clearly, we have
$$(6)  \quad  H^*(Hom_{kC_G(x)}(Bar_*(kG)\otimes_{kG}k,k))\simeq
H^*(Hom_{kC_G(x)}(Bar_*(kC_G(x))\otimes_{kC_G(x)}k,k)).$$ To give an
explicit isomorphism in $(6)$, we need to construct the comparison
maps between two projective resolutions $Bar_*(kG)\otimes_{kG}k$ and
$Bar_*(kC_G(x))\otimes_{kC_G(x)}k$ of the trivial $kC_G(x)$-module
$k$.

The comparison map from $Bar_*(kC_G(x))\otimes_{kC_G(x)}k$ to
$Bar_*(kG)\otimes_{kG}k$ is just the inclusion map
$$\iota: kC_G(x)\otimes \overline{kC_G(x)}^{\otimes
n}\hookrightarrow kG\otimes \overline{kG}^{\otimes n}.$$ This is
obvious or can be obtained using a setwise self-homotopy on
$Bar_*(kG)\otimes_{kG}k$ (see below for its explicit form).

 To
construct the comparison map on the reverse direction, we use a
setwise self-homotopy over $kC_G(x)\otimes
\overline{kC_G(x)}^{\otimes *}$ as follows (for $h_0\in C_G(x),
h_1,\cdots, h_{n}\in \overline{C_G(x)}$):
$$kC_G(x)\otimes \overline{kC_G(x)}^{\otimes
n}\longrightarrow kC_G(x)\otimes \overline{kC_G(x)}^{\otimes n+1},$$
$$h_0\otimes h_1\otimes \cdots \otimes h_{n}\longmapsto
1\otimes h_0\otimes h_1\otimes \cdots \otimes h_{n}.$$ Then we get
a comparison map $$\rho: Bar_*(kG)\otimes_{kG}k\longrightarrow
Bar_*(kC_G(x))\otimes_{kC_G(x)}k$$ as follows:
$$\rho_{-1}: k\longrightarrow k, \quad 1\longmapsto 1,$$
$$\rho_{0}: kG\longrightarrow kC_G(x), \quad h\gamma_{i,x}\longmapsto h,
\mbox{ for  }h\in  C_G(x),$$
$$\rho_{1}: kG\otimes \overline{kG}\longrightarrow kC_G(x)\otimes \overline{kC_G(x)}, \quad h\gamma_{i,x}\otimes g_1\longmapsto
h\otimes h_{i,1},$$
$$\mbox{ where
}\gamma_{i,x}g_1=h_{i,1}\gamma_{s_i,x}\mbox{ for }h_{i,1}\in
\overline{C_G(x)},$$
$$\cdots \cdots \cdots$$
$$\rho_{n}: kG\otimes \overline{kG}^{\otimes n}\longrightarrow kC_G(x)\otimes
 \overline{kC_G(x)}^{\otimes n}, \quad h\gamma_{i,x}\otimes g_1\otimes \cdots
\otimes g_n\longmapsto h\otimes h_{i,1}\otimes \cdots \otimes
h_{i,n},$$ where $ h_{i,1},\cdots, h_{i,n}\in \overline{C_G(x)}$
 are determined by the sequence $\{g_1,\cdots, g_n\}$ as follows:
$$\gamma_{i,x}g_1=h_{i,1}\gamma_{s_i^1,x}, \quad \gamma_{s_i^1,x}g_2=h_{i,2}
\gamma_{s_i^2,x}, \quad \cdots,
 \quad \gamma_{s_i^{n-1},x}g_n=h_{i,n}\gamma_{s_i^n,x}.$$ Notice that
$\rho\circ \iota=Id$ and $\iota\circ \rho\neq Id$.
 It follows
that we have two homomorphisms:
$$Hom_{kC_G(x)}(Bar_*(kG)\otimes_{kG}k,k)\longrightarrow
Hom_{kC_G(x)}(Bar_*(kC_G(x))\otimes_{kC_G(x)}k,k),$$
$$(\varphi_{5,x}: S_x\times\overline{G}^{\times n}\longrightarrow k)\longmapsto
(\varphi_{6,x}: \overline{C_G(x)}^{\times n}\longrightarrow k),
 \quad \varphi_{6,x}(h_1,\cdots, h_{n})= \varphi_{5,x}(1,h_1,\cdots,
h_{n})=a_{1,x},
$$
$$\mbox{ where }a_{1,x}\mbox{ is the cofficients of }x\mbox{ in }\varphi_{3,x}(h_1,\cdots,
h_n)=\sum_{i=1}^{n_x}a_{i,x}x_i;$$ and
$$Hom_{kC_G(x)}(Bar_*(kC_G(x))\otimes_{kC_G(x)}k,k)\longrightarrow
Hom_{kC_G(x)}(Bar_*(kG)\otimes_{kG}k,k),$$
$$(\varphi_{6,x}: \overline{C_G(x)}^{\times n}\longrightarrow k)\longmapsto
(\varphi_{5,x}: S_x\times\overline{G}^{\times n}\longrightarrow k),
 \quad \varphi_{5,x}(\gamma_{i,x},g_1,\cdots, g_{n})=\varphi_{6,x}(h_{i,1},
\cdots, h_{i,n}),
$$
 where for $ h_{i,1},\cdots, h_{i,n}\in \overline{C_G(x)}$
 are determined by the sequence $\{g_1,\cdots, g_n\}$ as follows:
$$\gamma_{i,x}g_1=h_{i,1}\gamma_{s_i^1,x}, \quad \gamma_{s_i^1,x}g_2=h_{i,2}\gamma_{s_i^2,x}, \quad \cdots, \quad
\gamma_{s_i^{n-1},x}g_n=h_{i,n}\gamma_{s_i^n,x}.$$ Since both $\iota$
and $\rho$ induce the identity map $1: k\longrightarrow k$, by Lemma
\ref{homotopy-equivalence}, we have inverse isomorphisms between
$H^*(Hom_{kC_G(x)}(Bar_*(kG)\otimes_{kG}k,k))$ and
$H^*(Hom_{kC_G(x)}(Bar_*(kC_G(x))\otimes_{kC_G(x)}k,k))$. The
correspondence is induced by $\varphi_{5,x}\longleftrightarrow
\varphi_{6,x}$, as we stated above. So we realize an isomorphism in
$(6)$ and its inverse.

\medskip
Summarizing the above six steps, we get the following main result in
this section.

\begin{Thm}\label{realization}  Let $k$ be a field and $G$ a finite group. Consider the additive
decomposition of Hochschild cohomology ring of the group algebra
$kG$:
$$HH^*(kG)\simeq \bigoplus_{x\in X}H^*(C_G(x),k),$$
where $X$ is a set of representatives of conjugacy classes of
elements of $G$ and $C_G(x)$ is the centralizer subgroup of $G$. We
compute the Hochschild cohomology $HH^*(kG)=H^*(Hom_{k(G\times
G)}(Bar_*(kG),kG))$ by the classical normalized bar resolution, and
we compute the group cohomology $H^*(C_G(x),k)$ by
$H^*(Hom_{kC_G(x)}(Bar_*(kC_G(x))\otimes_{kC_G(x)}k,k))$. Then, we
can realize an isomorphism in additive decomposition as follows:
 $$HH^*(kG)\stackrel{\sim}{\longrightarrow} \bigoplus_{x\in X}H^*(C_G(x),k),$$
$$[\varphi: \overline{G}^{\times n}\longrightarrow kG]\longmapsto [\widehat{\varphi}]=\bigoplus_{x\in X}[\widehat{\varphi}_x], \quad \widehat{\varphi}_x:
\overline{C_G(x)}^{\times n}\longrightarrow k,$$
$$\widehat{\varphi}_x(h_1,\cdots, h_{n})=a_{1,x}, \mbox{ where }
\pi_x((-1)^{\frac{n(n+1)}{2}}h_1\cdots h_n\varphi(h_{n}^{-1},\cdots,
h_1^{-1}))=\sum_{i=1}^{n_x}a_{i,x}x_i.$$ In other word,
$\widehat{\varphi}_x(h_1,\cdots, h_{n})$ is just the coefficient of
$x$ in $(-1)^{\frac{n(n+1)}{2}}h_1\cdots
h_n\varphi(h_{n}^{-1},\cdots, h_1^{-1})\in kG$. The inverse of the
above isomorphism is given as follows:
 $$\bigoplus_{x\in X}H^*(C_G(x),k)\stackrel{\sim}{\longrightarrow} HH^*(kG),$$
$$[\widehat{\varphi}]=\bigoplus_{x\in X}[\widehat{\varphi}_x], \quad \widehat{\varphi}_x:
\overline{C_G(x)}^{\times n}\longrightarrow k\longmapsto [\varphi:
\overline{G}^{\times n}\longrightarrow kG],$$
$$\varphi(g_1,\cdots, g_n)=(-1)^{\frac{n(n+1)}{2}}g_1\cdots g_n\sum_{x\in X}\sum_{i=1}^{n_x}\widehat{\varphi}_x(h_{i,1}',\cdots,
h_{i,n}')x_i,$$
$$\mbox{ where for }x\in X, h_{i,1}',\cdots, h_{i,n}'\in \overline{C_G(x)} \mbox{ are determined by the sequence }\{g_n^{-1},\cdots,
g_1^{-1}\}\mbox{ as follows: }$$
$$ \gamma_{i,x}g_n^{-1}=h_{i,1}'\gamma_{s_i^1,x},\quad \gamma_{s_i^1,x}g_{n-1}^{-1}=h_{i,2}'\gamma_{s_i^2,x},\quad \cdots, \quad \gamma_{s_i^{n-1},x}g_1^{-1}=h_{i,n}'\gamma_{s_i^n,x}.$$

\end{Thm}

\begin{Proof}
This is a direct consequence by applying the above isomorphisms
from $(1)$ to $(6)$ and their inverses. For an element $\varphi:
\overline{G}^{\times n}\longrightarrow kG$ in the $n$-th term
$C^n(kG)\simeq Map(\overline{G}^{\times n},kG)$ of the Hochschild
cohomology complex, $[\varphi]$ denotes the corresponding element
in the Hochschild cohomology group $HH^n(kG)$. Note that the
elements $h_{i,1}',\cdots, h_{i,n}'$ depend on $x\in X$ and the
sequence $\{g_n^{-1},\cdots,
g_1^{-1}\}$. For the simplicity of notations, we
avoid to write them down explicitly.
\end{Proof}

\begin{Rem}\label{illustration}
$(a)$ The correspondence in Theorem \ref{realization} makes use of
the same line employed  by Siegel and Witherspoon in
\cite{SW1999}. The difference is: they realize each step between
cohomology groups using standard operations like restriction,
induction, conjugation, etc., while we construct maps directly in
each step on the cohomology complex level.

\medskip
$(b)$  In \cite{SW1999}, as the authors  proved that
$HH^*(kG)\simeq H^*(G, {}_C kG)$ as graded algebras, they
concentrated on $H^*(G, {}_c kG)$ instead of $HH^*(kG)$ in most
part of their paper.  If we only consider the isomorphisms
(2)-(5), then
 the correspondence in Theorem
\ref{realization} become simpler:
$$H^*(G, {}_ckG)\stackrel{\sim}{\longrightarrow} \bigoplus_{x\in X}H^*(C_G(x),k),$$
$$[\varphi: \overline{G}^{\times n}\longrightarrow kG]\longmapsto [\widehat{\varphi}]=\bigoplus_{x\in X}[\widehat{\varphi}_x], \quad \widehat{\varphi}_x:
\overline{C_G(x)}^{\times n}\longrightarrow k,$$
$$\widehat{\varphi}_x(h_1,\cdots, h_{n})=a_{1,x}, \mbox{ the coefficient of
}x\mbox{ in } \varphi(h_1,\cdots, h_n)\in kG;$$
 $$\bigoplus_{x\in X}H^*(C_G(x),k)\stackrel{\sim}{\longrightarrow} H^*(G, {}_ckG),$$
$$[\widehat{\varphi}]=\bigoplus_{x\in X}[\widehat{\varphi}_x], \quad \widehat{\varphi}_x:
\overline{C_G(x)}^{\times n}\longrightarrow k\longmapsto [\varphi:
\overline{G}^{\times n}\longrightarrow kG],$$
$$\varphi(g_1,\cdots, g_n)=\sum_{x\in X}\sum_{i=1}^{n_x}\widehat{\varphi}_x(h_{i,1},\cdots,
h_{i,n})x_i,$$ where for $x\in X, h_{i,1},\cdots, h_{i,n}\in
\overline{C_G(x)}$ are determined by the sequence $\{g_1,\cdots,
g_n\}$ as follows:
$$\gamma_{i,x}g_1=h_{i,1}\gamma_{s_i^1,x}, \quad \gamma_{s_i^1,x}g_2=h_{i,2}\gamma_{s_i^2,x}, \quad \cdots, \quad
\gamma_{s_i^{n-1},x}g_n=h_{i,n}\gamma_{s_i^n,x}.$$
\end{Rem}

\bigskip


\section{Another realization of  the additive decomposition}

In \cite{CS1997}, Cibils and Solotar constructed a subcomplex of the
Hochschild cohomology complex for each conjugacy class, and then
they showed that for a finite abelian group, the subcomplex is
isomorphic to the complex computing group cohomology. We will
generalize this to any finite group: for each conjugacy class, this
complex computes the cohomology of the corresponding centralizer
subgroup. As a result, we give a second way to realize the additive
decomposition.

As before, let $k$ be a field and $G$ a finite group. Recall that
the Hochschild cohomology $HH^*(kG)$ of the group algebra $kG$ can
be computed by the following (cochain) complex:
$$(\mathcal{H}^*)\quad \quad \quad \quad 0\longrightarrow kG\stackrel{\delta_0}{\longrightarrow} Map(\overline{G},kG)\stackrel{\delta_1}{\longrightarrow} \cdots \longrightarrow Map(\overline{G}^{\times n},kG)
\stackrel{\delta_n}{\longrightarrow} \cdots,$$ where the differential is
given by $$\delta_0(x)(g)=gx-xg  \quad (\mbox{for } x\in kG \mbox{ and }
g\in \overline{G})$$ and (for $\varphi: \overline{G}^{\times
n}\longrightarrow kG$ and $g_1,\cdots, g_{n+1}\in \overline{G}$)
$$\quad  \quad \quad \delta_n(\varphi)(g_1,\cdots, g_{n+1})=g_1\varphi(g_2,\cdots,
g_{n+1}) + \quad  \quad  \quad  \quad  \quad  \quad  \quad  \quad
\quad  \quad  \quad  \quad  \quad  \quad  \quad  \quad  \quad  \quad
\quad  \quad \quad  \quad  \quad  \quad  \quad  \quad  \quad  \quad
\quad \quad  \quad $$
$$\quad  \quad  \quad  \quad  \quad  \quad  \quad  \quad
\quad  \quad \quad  \quad  \quad  \quad
\sum_{i=1}^n(-1)^i\varphi(g_1,\cdots, g_ig_{i+1},\cdots,
g_{n+1})+(-1)^{n+1}\varphi(g_1,\cdots, g_n)g_{n+1}.$$ We keep the
following notations in Section 3:  $X$ is a complete set of
representatives of the conjugacy classes in the finite group $G$.
For $x\in X$, $C_x=\{gxg^{-1}|g\in G\}$ is the conjugacy class
corresponding to $x$ and $C_G(x)=\{g\in G|gxg^{-1}=x\}$ is the
centralizer subgroup. Now take a conjugacy class $C_x$ and define
$$\mathcal{H}_x^0=kC_x,  \mbox{ and for }n\geq 1,$$
$$\mathcal{H}_x^n=\{\varphi: \overline{G}^{\times
n}\longrightarrow kG|\varphi(g_1,\cdots, g_{n})\in k[g_1\cdots
g_{n}C_x]\subset kG, \forall g_1,\cdots, g_n\in \overline{G}\},$$
where $g_1\cdots g_{n}C_x$ denotes the subset of $G$ by multiplying
$g_1\cdots g_{n}$ on $C_x$ and $k[g_1\cdots g_{n}C_x]$ is the
$k$-subspace of $kG$ generated by this set. Note that we have $g_1\cdots g_{n}C_x=C_xg_1\cdots g_{n}$
 and $k[g_1\cdots g_{n}C_x]=k[C_xg_1\cdots g_{n}]$.
Let
$\mathcal{H}_x^*=\bigoplus_{n\geq 0}\mathcal{H}_x^n$. Cibils and
Solotar (\cite[Page 20, Proof of the theorem]{CS1997}) observed that
$\mathcal{H}_x^*$ is a subcomplex of $\mathcal{H}^*$ and
$\mathcal{H}^*=\bigoplus_{x\in X}\mathcal{H}_x^*$.

\begin{Lem}\label{isomorphism}  $\mathcal{H}_x^*$ is canonically isomorphic to the
complex $Hom_{kG}(Bar_*(kG)\otimes_{kG}k,kC_x)$, which computes the
group cohomology $H^*(G,kC_x)$ of $G$ with coefficients in $kC_x$,
where $kC_x$ is a left $kG$-module under conjugation. \end{Lem}

\begin{Proof} We know from Section 3 that the complex
$Hom_{kG}(Bar_*(kG)\otimes_{kG}k,kC_x)$ is identified as the
following complex: $$0\longrightarrow
kC_x\stackrel{\delta_0}{\longrightarrow} Map(\overline{G},kC_x)
\stackrel{\delta_1}{\longrightarrow} \cdots \longrightarrow
Map(\overline{G}^{\times n},kC_x) \stackrel{\delta_n}{\longrightarrow}
\cdots,$$ where the differential is given by $$\delta_0(x)(g)=gxg^{-1}-x
\quad (\mbox{for } x\in kC_x \mbox{ and } g\in \overline{G})$$ and
(for $\varphi: \overline{G}^{\times n}\longrightarrow kC_x$ and
$g_1,\cdots, g_{n+1}\in \overline{G}$)
$$\quad  \quad \quad  \delta_n(\varphi)(g_1,\cdots, g_{n+1})=g_1\varphi(g_2,\cdots,
g_{n+1})g_1^{-1}+ \quad  \quad  \quad  \quad  \quad  \quad  \quad
\quad \quad  \quad  \quad  \quad  \quad  \quad  \quad  \quad  \quad
\quad \quad  \quad \quad  \quad  \quad  \quad  \quad  \quad  \quad
\quad \quad \quad  \quad $$ $$\quad  \quad  \quad  \quad  \quad
\quad \quad  \quad \quad  \quad \quad  \quad  \quad \quad
\sum_{i=1}^n(-1)^i\varphi(g_1,\cdots, g_ig_{i+1},\cdots,
g_{n+1})+(-1)^{n+1}\varphi(g_1,\cdots, g_n).$$
 A direct computation shows that the following map is an isomorphism of complexes:
$$\mathcal{H}_x^*\longrightarrow
Hom_{kG}(Bar_*(kG)\otimes_{kG}k,kC_x),$$
$$(\varphi_1: \overline{G}^{\times n}\longrightarrow kG)\longmapsto (\varphi_2:
\overline{G}^{\times n}\longrightarrow kC_x), \quad
\varphi_2(g_1,\cdots, g_{n})= \varphi_1(g_1,\cdots,
g_{n})g_n^{-1}\cdots g_1^{-1}.$$ Its inverse is given by
$$Hom_{kG}(Bar_*(kG)\otimes_{kG}k,kC_x)\longrightarrow
\mathcal{H}_x^*,$$
$$(\varphi_2: \overline{G}^{\times n}\longrightarrow kC_x)\longmapsto (\varphi_1:
\overline{G}^{\times n}\longrightarrow kG), \quad
\varphi_1(g_1,\cdots, g_{n})= \varphi_2(g_1,\cdots, g_{n})g_1\cdots
g_{n}.$$ Passing to the cohomology, we have
$H^*(\mathcal{H}_x^*)\simeq H^*(G,kC_x)$.
\end{Proof}

\begin{Rem} \label{ComparisonRealization} Since the first three steps of the previous section
realize
$$HH^*(kG)\simeq H^*(G, {}_ckG)\simeq \oplus_{x\in X} H^*(C_G(x),
k),$$  these isomorphisms also give  a decomposition of the
complex $\mathcal{H}^*$, which computes $HH^*(kG)$. In fact,
during these three steps, we establish the following isomorphisms
of complexes
$$\begin{array}{rcl}\mathcal{H}^*&=&Hom_{k(G\times G}(Bar_*(kG), kG)\\
&\stackrel{(1)}{\simeq}& Hom_{k(G\times G)}(k(G\times
G)\ot_{kG}Bar_*(kG)\ot_{kG}k, kG)\\
&\stackrel{(2)}{\simeq}& Hom_{kG}( Bar_*(kG)\ot_{kG}k,
{}_ckG)\\
&\stackrel{(3)}{=}&\bigoplus_{x\in X}Hom_{kG}( Bar_*(kG)\ot_{kG}k,
{}_ckC_x).\end{array}$$ So the complex $Hom_{kG}(
Bar_*(kG)\ot_{kG}k, {}_ckC_x)$ is isomorphic to a subcomplex of
$\mathcal{H}^*$ and we verify easily that this subcomplex is just
the above defined $\mathcal{H}^*_x$. However, the isomorphism
between these two complexes  is as follows:
$$\mathcal{H}_x^*\longrightarrow
Hom_{kG}(Bar_*(kG)\otimes_{kG}k,kC_x),$$
$$(\varphi_1: \overline{G}^{\times n}\longrightarrow kG)\longmapsto (\varphi_2:
\overline{G}^{\times n}\longrightarrow kC_x),
 \quad \varphi_2(g_1,\cdots, g_{n})= (-1)^{\frac{n(n+1)}{2}}g_1\cdots
g_n\varphi_1(g_n^{-1},\cdots, g_1^{-1} ).$$ Its inverse is given by
$$Hom_{kG}(Bar_*(kG)\otimes_{kG}k,kC_x)\longrightarrow
\mathcal{H}_x^*,$$
$$(\varphi_2: \overline{G}^{\times n}\longrightarrow kC_x)\longmapsto (\varphi_1:
\overline{G}^{\times n}\longrightarrow kG), \quad
\varphi_1(g_1,\cdots, g_{n})= (-1)^{\frac{n(n+1)}{2}}g_1\cdots
g_n\varphi_1(g_n^{-1},\cdots, g_1^{-1} ).$$ Note that this
isomorphism differs from the one in Lemma~\ref{isomorphism} by an
automorphism of the complex $Hom_{kG}( Bar_*(kG)\ot_{kG}k,
{}_ckC_x)$, which sends $\varphi: \ol{G}^{\times n}\to kC_x$ to
$\varphi': \ol{G}^{\times n}\to kC_x$ with
$$\varphi'(g_1, \cdots, g_n)= (-1)^{\frac{n(n+1)}{2}}g_1\cdots
g_n\varphi(g_n^{-1},\cdots, g_1^{-1} ) g_n^{-1}\cdots
g_1^{-1}.$$

\end{Rem}

\medskip

On the other hand, we have shown that the complex
$Hom_{kG}(Bar_*(kG)\otimes_{kG}k,kC_x)$ is isomorphic to the complex
$Hom_{kC_G(x)}(Bar_*(kC_G(x))\otimes_{kC_G(x)}k,k)$, which computes
the group cohomology $H^*(C_G(x),k)$ of the centralizer subgroup
$C_G(x)$ with coefficients in the trivial module $k$. (cf. Section
3, from the fourth step to the six step.) Therefore we get another
realization to the additive decomposition:

\begin{Thm}\label{second-realization}  Let $k$ be a field and $G$ a finite group. Consider the additive
decomposition of Hochschild cohomology ring of the group algebra
$kG$:
$$HH^*(kG)\simeq \bigoplus_{x\in X}H^*(C_G(x),k)$$
where $X$ is a set of representatives of conjugacy classes of
elements of $G$ and $C_G(x)$ is the centralizer subgroup of $G$. We
compute the Hochschild cohomology $HH^*(kG)=H^*(Hom_{k(G\times
G)}(Bar_*(kG),kG))$ by the classical normalized bar resolution, and
we compute the group cohomology $H^*(C_G(x),k)$ by
$H^*(Hom_{kC_G(x)}(Bar_*(kC_G(x))\otimes_{kC_G(x)}k,k))$. Then, we
can realize an isomorphism in additive decomposition as follows:
 $$HH^*(kG)\stackrel{\sim}{\longrightarrow} \bigoplus_{x\in X}H^*(C_G(x),k),$$
$$[\varphi_x: \overline{G}^{\times n}\longrightarrow kG], \quad \varphi_x\in\mathcal{H}_x^n\longmapsto [\widehat{\varphi}_x:
\overline{C_G(x)}^{\times n}\longrightarrow k],$$
$$\widehat{\varphi}_x(h_1,\cdots, h_{n})=a_{1,x}, \mbox{ where }
\varphi_x(h_1,\cdots, h_n)h_{n}^{-1}\cdots
h_1^{-1}=\sum_{i=1}^{n_x}a_{i,x}x_i\in kC_x.$$ In other word,
$\widehat{\varphi}_x(h_1,\cdots, h_{n})$ is just the coefficient of
$x$ in $\varphi_x(h_1,\cdots, h_n)h_{n}^{-1}\cdots h_1^{-1}\in
kC_x$. The inverse of the above isomorphism is given as follows:
 $$\bigoplus_{x\in X}H^*(C_G(x),k)\stackrel{\sim}{\longrightarrow} HH^*(kG),$$
$$[\widehat{\varphi}_x:
\overline{C_G(x)}^{\times n}\longrightarrow k]\longmapsto
[\varphi_x: \overline{G}^{\times n}\longrightarrow kG], \quad
\varphi_x\in\mathcal{H}_x^n,$$
$$\varphi_x(g_1,\cdots, g_n)=\sum_{i=1}^{n_x}\widehat{\varphi}_x(h_{i,1},\cdots,
h_{i,n})x_ig_1\cdots g_n,$$
$$\mbox{ where }h_{i,1},\cdots, h_{i,n}\in \overline{C_G(x)} \mbox{ are determined by the sequence }\{g_1,\cdots,
g_n\}\mbox{ as follows: }$$
$$\gamma_{i,x}g_1=h_{i,1}\gamma_{s_i^1,x}, \quad \gamma_{s_i^1,x}g_2=h_{i,2}\gamma_{s_i^2,x}, \quad \cdots, \quad
 \gamma_{s_i^{n-1},x}g_n=h_{i,n}\gamma_{s_i^n,x}.$$

\end{Thm}

\begin{Proof}
This is a combination of Lemma \ref{isomorphism} and the
correspondence from the fourth step to the six step in Section 3.
\end{Proof}

  By Remark~\ref{ComparisonRealization},  the two realizations of the additive
decomposition  in   Theorem~\ref{realization} and Theorem
\ref{second-realization} are essentially the same on the cohomology group
level. In the sequel, we prefer to the second realization since it
is simpler.

\bigskip

\section{The cup product formula}

 We keep the notations of the previous sections: $k$ is a field, and $G$ is a finite group, and so on.
 We describe the cup product formula for the Hochschild cohomology ring
 $HH^*(kG)$ in terms of the additive decomposition.

 We shall define a product over $\sum_{x\in X}H^*(C_G(x), k)$ such
 that the isomorphism $$HH^*(kG)\simeq \bigoplus_{x\in X}H^*(C_G(x),
 k)$$ realized in Theorem~\ref{second-realization} becomes an isomorphism of graded algebras.

 Let
$[\widehat{\varphi_x}]\in H^*(C_G(x), k)$ (respectively
$[\widehat{\varphi_y}]\in H^m(C_G(y), k)$)  represented by the map
$\widehat{\varphi_x}: \overline{C_G(x)}^{\times n}\longrightarrow
k$ (respectively by $\widehat{\varphi_y}:
\overline{C_G(y)}^{\times m}\longrightarrow k$). Define
$$[ \widehat{\varphi_x}\cup  \widehat{\varphi_y}]=\sum_{z\in
X}[( \widehat{\varphi_x}\cup  \widehat{\varphi_y})_z]\in
\bigoplus_{z\in X}H^*(C_G(z),
 k)$$
 with $(\widehat{\varphi_x}\cup  \widehat{\varphi_y})_z:
 \ol{C_G(z)}^{\times n+m}\to k$  as follows:
$$(\widehat{\varphi_x}\cup \widehat{\varphi_y})_z(h_1,\cdots, h_{n},h_{n+1}, \cdots, h_{n+m})
=\sum_{(i,j)\in I_1'}\widehat{\varphi_x}(h_{i,1}, \cdots,
h_{i,n})\widehat{\varphi_y}(h_{j,1}, \cdots, h_{j,m}),$$ where
\begin{itemize}\item $I_1'$ is the set of pairs $(i, j)$ such that
$$x_ih_1\cdots h_n y_j (h_{1}\cdots h_{n})^{-1}=z;$$

\item $$\gamma_{i, x}h_1=h_{i,1} \gamma_{s_i^1, x}, \quad \gamma_{s_i^1,
x}h_2=h_{i,2} \gamma_{s_i^2, x}, \quad \cdots, \quad
\gamma_{s_i^{n-1}, x}h_n=h_{i,n} \gamma_{s_i^n, x};$$

\item $$\gamma_{j, y}h_{n+1}=h_{j,1} \gamma_{s_j^1, y}, \quad
\gamma_{s_j^1, y}h_{n+2}=h_{j,2} \gamma_{s_j^2, y}, \quad \cdots,
 \quad \gamma_{s_j^{m-1}, j}h_{n+m}=h_{j,m} \gamma_{s_j^m, y}.$$

\end{itemize}


\begin{Thm}\label{cup-product} With the product defined above, the isomorphism  in Theorem
\ref{second-realization} is an isomorphism of graded algebras.
\end{Thm}

\begin{Proof}

We shall show that with respect to the isomorphism in Theorem
\ref{second-realization}, the product defined above coincide with
the cup product on the cohomology complex level.

 Let $[\widehat{\varphi_x}]
(\widehat{\varphi_x}: \overline{C_G(x)}^{\times n}\longrightarrow
k)$ and $[\widehat{\varphi_y}] (\widehat{\varphi_y}:
\overline{C_G(y)}^{\times m}\longrightarrow k)$ be two elements in
$H^n(C_G(x),k)$ and in $H^m(C_G(y),k)$, respectively. Denote by
$[\varphi_x: \overline{G}^{\times n}\longrightarrow kG]
(\varphi_x\in\mathcal{H}_x^n)$ and $[\varphi_y:
\overline{G}^{\times m}\longrightarrow kG]
(\varphi_y\in\mathcal{H}_y^m)$ be the corresponding elements in
$HH^*(kG)$. By Theorem~\ref{second-realization},
$$\begin{array}{rll}\varphi_x&:&
\overline{G}^{\times n}\longrightarrow kG,\\
 (g_1, \cdots, g_n)&\mapsto&  \sum_{i=1}^{n_x}  \widehat{\varphi_x}(h_{i,1},
\cdots, h_{i,n})x_i g_1\cdots g_n,\end{array}$$ with
$$\gamma_{i, x}g_1=h_{i,1} \gamma_{s_i^1, x}, \quad \gamma_{s_i^1,
x}g_2=h_{i,2} \gamma_{s_i^2, x}, \quad \cdots, \quad
\gamma_{s_i^{n-1}, x}g_n=h_{i,n} \gamma_{s_i^n, x}. $$  A  similar
formula  for $\varphi_y$ works as well.

 Now  denote by $\varphi_x\cup
 \varphi_y: \overline{G}^{\times (n+m)}\longrightarrow
kG$ the cup product. By the definition of the cup product,  for
any $z\in X$, we obtain $(\varphi_x\cup
 \varphi_y)_z\in\mathcal{H}_z^{n+m}$ given by
$$\begin{array}{rll}( \varphi_x\cup
 \varphi_y)_z&:& \overline{G}^{\times
(n+m)}\longrightarrow kG, \\
 (g_1,\cdots, g_{n},\cdots, g_{n+m})&\mapsto& \sum_{k=1}^{n_z}\sum_{(i,j)\in I_k} \widehat{\varphi_x}(h_{i,1},\cdots,
h_{i,n}) \widehat{\varphi_y}(h_{j,1},\cdots, h_{j,m})z_k g_1\cdots
g_{n+m},\end{array}$$ where \begin{itemize}\item $I_k$ is the set
of pairs $(i, j)$ such that
$$x_ig_1\cdots g_n y_j (g_{1}\cdots g_{n})^{-1}=z_k;$$

\item $$\gamma_{i, x}g_1=h_{i,1} \gamma_{s_i^1, x}, \quad \gamma_{s_i^1,
x}g_2=h_{i,2} \gamma_{s_i^2, x}, \quad \cdots, \quad
\gamma_{s_i^{n-1}, x}g_n=h_{i,n} \gamma_{s_i^n, x};$$

\item $$\gamma_{j, y}g_{n+1}=h_{j,1} \gamma_{s_j^1, y}, \quad
\gamma_{s_j^1, y}g_{n+2}=h_{j,2} \gamma_{s_j^2, y}, \quad \cdots,
 \quad \gamma_{s_j^{m-1}, j}g_{n+m}=h_{j,m} \gamma_{s_j^m, y}.$$

\end{itemize}
Note that $I_k$ depends on the elements $g_1,\cdots,g_n$.
Again by Theorem~\ref{second-realization}, we obtain an element in
$H^*(C_G(z), k)$   of the following form:
$$ \begin{array}{rll}\overline{C_G(z)}^{\times
(n+m)}&\longrightarrow& k\\
 (h_1,\cdots, h_{n},h_{n+1}, \cdots, h_{n+m})&\longmapsto&  \sum_{(i,j)\in I_1}\widehat{\varphi_x}(h_{i, 1}, \cdots,
h_{i,n})\widehat{\varphi_y}(h_{j,1}, \cdots,
h_{j,m}),\end{array}$$ which is just $(\widehat{\varphi_x} \cup
\widehat{\varphi_y})_z$ defined before.

\end{Proof}

Similarly we can prove the following result.

\begin{Thm}\label{cup-product-2} The isomorphism  in Theorem
\ref{realization} is an isomorphism of graded algebras with
respect to the  following product defined on  $\bigoplus_{x\in
X}H^m(C_G(x),k)$.

 Let
$[\widehat{\varphi}_x] (\widehat{\varphi}_x:
\overline{C_G(x)}^{\times n}\longrightarrow k)$ and
$[\widehat{\varphi}_y] (\widehat{\varphi}_y:
\overline{C_G(y)}^{\times m}\longrightarrow k)$ be two elements in
$H^n(C_G(x),k)$ and in $H^m(C_G(y),k)$, respectively. Define
$$[ \widehat{\varphi_x}\cup  \widehat{\varphi_y}]=\sum_{z\in
X}[(  \widehat{\varphi_x}\cup  \widehat{\varphi_y})_z]\in
\bigoplus_{z\in X}H^*(C_G(z),
 k)$$
 with $( \widehat{\varphi_x}\cup  \widehat{\varphi_y})_z:
 \ol{C_G(z)}^{\times n+m}\longrightarrow k$  as follows:
$$( \widehat{\varphi_x}\cup  \widehat{\varphi_y})_z(h_1,\cdots, h_{n},h_{n+1}, \cdots, h_{n+m})
=(-1)^{nm}\sum_{(i,j)\in I_1}\widehat{\varphi_x}(h_{i,1}', \cdots,
h_{i,n}')\widehat{\varphi_y}(h_{j,1}', \cdots, h_{j,m}'),$$ where
\begin{itemize}\item $I_1$ is the set of pairs $(i, j)$ such that
$$ h_1\cdots h_mx_i  (h_1\cdots h_m)^{-1}y_j=z_1=z;$$

\item $$\gamma_{i, x}h_{m+1}=h_{i,1}' \gamma_{s_i^1, x}, \quad
\gamma_{s_i^1, x}h_{m+2}=h_{i,2}' \gamma_{s_i^2, x}, \quad \cdots,
\quad \gamma_{s_i^{n-1}, x}h_{n+m}=h_{i,n}' \gamma_{s_i^n, x};$$

\item $$\gamma_{j, y}h_{1}=h_{j,1}' \gamma_{s_j^1, y}, \quad
\gamma_{s_j^1, y}h_{2}=h_{j,2}' \gamma_{s_j^2, y}, \quad  \cdots,
\quad \gamma_{s_j^{m-1}, j}h_{m}=h_{j,m}' \gamma_{s_j^m, y}.$$

\end{itemize}

\end{Thm}

\begin{Rem} (1) By Remark \ref{illustration}
$(a)$, our cup product formulae in Theorems \ref{cup-product} and
\ref{cup-product-2} are consistent with Siegel and Witherspoon's
formula in \cite[Theorem 5.1]{SW1999} up to an isomorphism.

(2) From our realization of the graded algebra isomorphism
$$HH^*(kG)\simeq \bigoplus_{x\in X}H^*(C_G(x), k)=H^*(G,
k)\oplus(\bigoplus_{x\in X- \{1\}}H^*(C_G(x), k)),$$ it is clear
that $H^*(G, k)$ can be seen as a graded subalgebra of $HH^*(kG)$
and each $H^*(C_G(x), k)$ is a graded $H^*(G, k)$-submodule of
$HH^*(kG)$. Therefore, the additive decomposition gives an
isomorphism of graded $H^*(G, k)$-modules.
\end{Rem}

\bigskip

\section{The $\bigtriangleup$ operator formula}

Let $k$ be a field and $G$ a finite group. Recall that the group
algebra $kG$ is a symmetric algebra with the bilinear form
$$\langle~,~\rangle: kG\times kG\longrightarrow k,$$
$$\langle g,h\rangle=\left\{\begin{array}{cc}
1 & \mbox{ if }g=h^{-1} \\ 0 & \mbox{ otherwise }
\end{array}\right. $$ for $g,h\in G$.
For $n\geq 1$, the operator $\bigtriangleup: HH^n(kG)\longrightarrow
HH^{n-1}(kG)$ on the Hochschild cohomology is defined by the
equation
$$\langle\bigtriangleup(\varphi)(g_1,
 \cdots,g_{n-1}),g_n\rangle=\sum_{i=1}^n(-1)^{i(n-1)}\langle
 \varphi(g_i,
 \cdots, g_{n-1}, g_n, g_1,
 \cdots, g_{i-1}),1\rangle,$$ where $\varphi\in C^n(kG)\simeq Map(\overline{G}^{\times n},kG)$,
$\bigtriangleup(\varphi)\in C^{n-1}(kG)\simeq
Map(\overline{G}^{\times n-1},kG)$. Equivalently,
$$ \bigtriangleup(\varphi)(g_1,
 \cdots,g_{n-1})=\sum_{g_n\in G}\sum_{i=1}^n(-1)^{i(n-1)}\langle \varphi(g_i, \cdots, g_{n-1}, g_n,
 g_1, \cdots, g_{i-1}),1\rangle g_n^{-1}.$$ This operator together with
the cup product $\cup$ and the Lie bracket $[~,~]$ defines a BV
algebra structure on $HH^*(kG)$.

We know from Section 4 that, for a conjugacy class $C_x$ of $G$,
$\mathcal{H}_x^*=\bigoplus_{n\geq 0}\mathcal{H}_x^n$ is a subcomplex
of the Hochschild cohomology complex $\mathcal{H}^*$, where
$$\mathcal{H}_x^n=\{\varphi: \overline{G}^{\times
n}\longrightarrow kG|\varphi(g_1,\cdots, g_{n})\in k[g_1\cdots
g_{n}C_x]\subset kG, \forall g_1,\cdots, g_n\in \overline{G}\}.$$

\begin{Lem}\label{restriction}  The operator $\bigtriangleup: \mathcal{H}^n\longrightarrow
\mathcal{H}^{n-1}$ restricts to $\bigtriangleup_x:
\mathcal{H}_x^n\longrightarrow \mathcal{H}_x^{n-1}$ for each
conjugacy class $C_x$.
\end{Lem}

\begin{Proof} We need to show that $\bigtriangleup(\varphi)\in \mathcal{H}_x^{n-1}$ for each $\varphi\in \mathcal{H}_x^n$.
Let  $g_1,\cdots, g_{n-1}\in \overline{G}$. Since $$
\bigtriangleup(\varphi)(g_1,
 \cdots,g_{n-1})=\sum_{g_n\in G}\sum_{i=1}^n(-1)^{i(n-1)}\langle \varphi(g_i, \cdots, g_{n-1}, g_n,
 g_1, \cdots, g_{i-1}),1\rangle g_n^{-1},$$
 it suffices to prove the
 following statement: if $\langle \varphi(g_i,
\cdots, g_{n-1}, g_n,
 g_1, \cdots, g_{i-1}),1\rangle\neq 0$ for some $i$, then $g_n^{-1}\in g_1\cdots
 g_{n-1}C_x$. Indeed, $\langle \varphi(g_i,
\cdots, g_{n-1}, g_n,
 g_1, \cdots, g_{i-1}),1\rangle\neq 0$ implies that $1\in g_i\cdots g_{n-1}g_ng_1\cdots
 g_{i-1}C_x$, or equivalently, $g_n^{-1}\in g_1\cdots g_{i-1}C_xg_i\cdots g_{n-1}=g_1\cdots
 g_{n-1}C_x$.
\end{Proof}

Now we can determine the behavior of the operator $\bigtriangleup$
under the additive decomposition.

\begin{Thm} \label{bigtriangleup-operator}
  Let
$\widehat{\bigtriangleup}_x: H^n(C_G(x),k)\longrightarrow
H^{n-1}(C_G(x),k)$ be the map induced by the operator
$\bigtriangleup_x: HH^n(kG)\longrightarrow HH^{n-1}(kG)$ via the
isomorphism established in Lemma~\ref{isomorphism}. Then
$\widehat{\bigtriangleup}_x$ is defined as follows:
$$\widehat{\bigtriangleup}_x(\psi)(h_1,\cdots, h_{n-1})=\sum_{i=1}^n(-1)^{i(n-1)}\psi(h_i, \cdots, h_{n-1}, h_{n-1}^{-1}\cdots h_1^{-1}x^{-1},
 h_1, \cdots, h_{i-1})$$
 for $\psi:\overline{C_G(x)}^{\times
n}\longrightarrow k$ and for $h_1,\cdots, h_{n-1}\in
\overline{C_G(x)}$.

\end{Thm}

\begin{Proof}We shall prove that  the following  diagram
$$\xymatrix{
H^n(\mathcal{H}_x^{*})\ar[r]^{\bigtriangleup_x}\ar[d]_{\wr}  &
H^{n-1}(\mathcal{H}_x^*)
\ar[d]_{\wr}\\
  H^n(C_G(x),k)\ar[r]^{\widehat{\bigtriangleup}_x}  &
H^{n-1}(C_G(x),k)}$$ is commutative, where the vertical
isomorphisms are given in Lemma~\ref{isomorphism}.

Take an element $\psi:\overline{C_G(x)}^{\times n}\longrightarrow
k$ in $Hom_{kC_G(x)}(Bar_n(kC_G(x))\otimes_{kC_G(x)}k,k)$ and
denote by $\varphi: \overline{G}^{\times n}\longrightarrow kG$ the
corresponding element in $\mathcal{H}_x^{n}$. By Theorem
\ref{second-realization}, for any $h_1,\cdots, h_{n}\in
\overline{C_G(x)}$, $\psi(h_1,\cdots, h_{n})$ is equal to the
coefficient of $x$ in $\varphi(h_1,\cdots, h_n)h_{n}^{-1}\cdots
h_1^{-1}\in kC_x$.  We should prove that
$\widehat{\bigtriangleup}_x(\psi)$ corresponds to
$\bigtriangleup_x(\varphi)$ via the isomorphism in
Lemma~\ref{isomorphism}.

Now
$$\bigtriangleup_x(\varphi)(g_1,
 \cdots,g_{n-1})=\sum_{g_n\in G}\sum_{i=1}^n(-1)^{i(n-1)}\langle \varphi(g_i, \cdots, g_{n-1}, g_n,
 g_1, \cdots, g_{i-1}),1\rangle g_n^{-1}.$$
   For any $h_1,\cdots, h_{n-1}\in
\overline{C_G(x)}$, consider  the coefficient of $x$ in
$\bigtriangleup_x(\varphi)(h_1,\cdots, h_{n-1})h_{n-1}^{-1}\cdots
h_1^{-1}\in kC_x$, or equivalently, the coefficient of $xh_1\cdots
h_{n-1}$ in $\bigtriangleup_x(\varphi)(h_1,\cdots, h_{n-1})\in
k[h_1\cdots h_{n-1}C_x]$. This coefficient is equal to
$$\langle\bigtriangleup_x(\varphi)(h_1,
 \cdots,h_{n-1}),h_{n-1}^{-1}\cdots
h_1^{-1}x^{-1}\rangle$$
$$\quad  \quad  \quad  \quad  \quad  \quad  \quad \quad  \quad  \quad
=\sum_{i=1}^n(-1)^{i(n-1)}\langle \varphi(h_i, \cdots, h_{n-1},
h_{n-1}^{-1}\cdots h_1^{-1}x^{-1},
 h_1, \cdots, h_{i-1}),1\rangle.$$

On the other hand, we also know that $\psi(h_i, \cdots, h_{n-1},
h_{n-1}^{-1}\cdots h_1^{-1}x^{-1},
 h_1, \cdots, h_{i-1})$ is equal to the
coefficient of $x$ in $$\varphi(h_i, \cdots, h_{n-1},
h_{n-1}^{-1}\cdots h_1^{-1}x^{-1},
 h_1, \cdots, h_{i-1})h_{i-1}^{-1}\cdots h_1^{-1}xh_1\cdots h_{n-1}h_{n-1}^{-1}\cdots h_i^{-1}$$
 $$=\varphi(h_i, \cdots, h_{n-1},
h_{n-1}^{-1}\cdots h_1^{-1}x^{-1},
 h_1, \cdots, h_{i-1})x\in kC_x,\quad  \quad \quad  \quad  \quad  \quad  \quad  \quad \quad  \quad  \quad \quad$$
which is again equal to $\langle \varphi(h_i, \cdots, h_{n-1},
h_{n-1}^{-1}\cdots h_1^{-1}x^{-1},
 h_1, \cdots, h_{i-1}),1\rangle$. It follows that
$$\widehat{\bigtriangleup}_x(\psi)(h_1,\cdots, h_{n-1})=\sum_{i=1}^n(-1)^{i(n-1)}\psi(h_i, \cdots, h_{n-1}, h_{n-1}^{-1}\cdots h_1^{-1}x^{-1},
 h_1, \cdots, h_{i-1}).$$

 We have proved that
$\widehat{\bigtriangleup}_x(\psi)$ corresponds to
$\bigtriangleup_x(\varphi)$ via the isomorphism in
Lemma~\ref{isomorphism} and the diagram is commutative (even at the
cohomology complex level).
\end{Proof}

\begin{Rem} \label{Definition-BV-sub-algebra} By \cite[Corollary 2.2]{FS2004}, we know that $H^*(G,k)$ is a
Gerstenhaber subalgebra of $HH^*(kG)$ under the inclusion map:
$$Hom_{kG}(Bar_*(kG)\otimes_{kG}k,k)\hookrightarrow Hom_{k(G\times
G)}(Bar_*(kG),kG),$$
$$(\varphi: \overline{G}^{\times n}\longrightarrow k)\longmapsto (\psi: \overline{G}^{\times n}\longrightarrow kG), \quad \psi(g_1,
\cdots, g_{n})=\varphi(g_1, \cdots, g_{n})g_1\cdots g_{n},$$ which
is in fact induced by the isomorphism in Lemma~\ref{isomorphism}
corresponding to $x=1$.

Notice that by notations in Section 6, $\psi\in \mathcal{H}_1^n$. So
motivated by Theorem \ref{bigtriangleup-operator}, we can similarly
define an operator $\bigtriangleup_1: H^n(G,k)\longrightarrow
H^{n-1}(G,k)$ in the group cohomology $H^*(G,k)$ as follows:
$$\bigtriangleup_1(\varphi)(g_1,\cdots, g_{n-1})=\sum_{i=1}^n(-1)^{i(n-1)}\varphi(g_i, \cdots, g_{n-1}, g_{n-1}^{-1}\cdots g_1^{-1},
 g_1, \cdots, g_{i-1})$$
 for $\varphi:\overline{G}^{\times
n}\longrightarrow k$ and for $g_1,\cdots, g_{n-1}\in \overline{G}$.

\end{Rem}

We prove that $H^*(G,k)$ is in fact a BV subalgebra of $HH^*(kG)$.

\begin{Cor} \label{BV-sub-algebra}  Let $k$ be a field and $G$ a finite group. Then $H^*(G,k)\hookrightarrow HH^*(kG)$ is a
BV subalgebra.

\end{Cor}

\begin{Proof} The inclusion in known to preserve the product
structure. In fact this is a direct consequence of
Theorem~\ref{cup-product}. In that result, taking $x=1=y$, we always
has $z=1$ if the set $I_1'$ is not empty.

As the inclusion $H^n(G,k)\hookrightarrow HH^*(kG)$ induced by the
isomorphism in Lemma~\ref{isomorphism} corresponding to $x=1$,
Theorem~\ref{bigtriangleup-operator} shows that it preserves the
$\bigtriangleup$-operator.  Since this operator together with the
cup product $\cup$ and the Lie bracket $[~,~]$ define a BV algebra
structure on $HH^*(kG)$,  via the isomorphism in
Lemma~\ref{isomorphism}, we deduce that the Lie bracket $[~,~]$
restricts to $H^*(G, k)=H^*(C_G(1), k)$.


\end{Proof}

Now we specialize to the case of abelian groups. Let $G$ be an
abelian group. In this case, the Hochschild cohomology ring
$HH^*(kG)$ of the group algebra $kG$ is isomorphic to the tensor
product algebra of $kG$ and the group cohomology ring $H^*(G,k)$:
$HH^*(kG)\simeq kG\otimes_kH^*(G,k).$   According to
\cite{CS1997}, this isomorphism is given as follows. For $G$ an
abelian group, conjugacy classes are elements of $G$, hence a
cochain $\varphi_x$ of $\mathcal{H}_x^n$ for $x\in G$ attributes a
scalar multiple of $g_1\cdots g_{n}x$ for each $(g_1,\cdots,
g_{n})\in \overline{G}^{\times n}$  and we obtain in this way a
map $\overline{\varphi_x}:\overline{G}^{\times n}\longrightarrow
k$. It is easy to see that
  the map $\widehat{\varphi_x}$ in
Theorem~\ref{second-realization} is just this scalar.

Now Theorem~\ref{cup-product} shows that  the map
$\varphi\longmapsto \Sigma_{x\in G} (x\otimes
\overline{\varphi_x})$ defines a ring isomorphism
$C^*(kG)\longrightarrow kG\otimes C^*(kG,k)$ compatible with the
differentials, and therefore it induces the above isomorphism.
Theorem~\ref{bigtriangleup-operator} specializes to the following
statement.

\begin{Prop} \label{abelian-bigtriangleup-operator} Let $k$ be a field and $G$ a finite abelian group.
Under the above isomorphism $HH^*(kG)\simeq kG\otimes_kH^*(G,k)$,
the operator $\bigtriangleup: HH^n(kG)\longrightarrow
HH^{n-1}(kG)$ corresponds to the sum of operators
$x\otimes\overline{\bigtriangleup_x}: x\otimes
H^n(G,k)\longrightarrow x\otimes H^{n-1}(G,k)$, where $x\in G$ and
$\overline{\bigtriangleup_x}: H^n(G,k)\longrightarrow
H^{n-1}(G,k)$ is defined as follows:
$$\overline{\bigtriangleup_x}(\overline{\varphi})(g_1,\cdots, g_{n-1})=
\sum_{i=1}^n(-1)^{i(n-1)}\overline{\varphi}(g_i, \cdots, g_{n-1}, g_{n-1}^{-1}\cdots g_1^{-1}x^{-1},
 g_1, \cdots, g_{i-1})$$
 for $\overline{\varphi}:\overline{G}^{\times
n}\longrightarrow k$ and for $g_1,\cdots, g_{n-1}\in \overline{G}$.

\end{Prop}

\begin{Rem} We could also use the first realization to deduce a
formula of the $\bigtriangleup$ operator. However, this formula is
much more complicated than that of
Theorem~\ref{bigtriangleup-operator}. We refrain from giving it
here.

\end{Rem}

In a BV-algebra, we have the following equation (see
\cite{Getzler1994}; Here we have changed the original equation
according to the sign convention in Remark \ref{sign-convention} and
we omit the sign $\cup$ in the equation):
$$\bigtriangleup(\alpha\beta\gamma)
=(-1)^{|\alpha||\beta||\gamma|}[(-1)^{|\gamma|}\bigtriangleup(\alpha\beta)\gamma
+ \alpha\bigtriangleup(\beta\gamma) +
(-1)^{|\alpha||\beta|}\beta\bigtriangleup(\alpha\gamma)$$
$$ - (-1)^{|\alpha|}
\bigtriangleup(\alpha)\beta\gamma -
(-1)^{|\alpha|+|\beta|-|\alpha||\gamma|}\alpha(\bigtriangleup(\beta))\gamma
-
(-1)^{|\alpha|+|\beta|+|\gamma|}\alpha\beta\bigtriangleup(\gamma)],$$
where $\alpha, \beta, \gamma\in HH^*(A)$ are homogeneous elements.
So in order to compute the $\bigtriangleup$ operator in $HH^*(A)$,
it suffices to find the value of $\bigtriangleup$ on each generator
and on the cup product of every two generators. Finally, let us
mention that we can use the cup product formula, the
$\bigtriangleup$ operator formula and the following formula to
compute the Lie bracket:
$$[\alpha,\beta]=-(-1)^{(|\alpha|-1)|\beta|}(\bigtriangleup(\alpha\cup\beta)-\bigtriangleup(\alpha)\cup\beta-(-1)^{|\alpha|}\alpha\cup\bigtriangleup(\beta)).$$

\bigskip

\section{The symmetric group of degree $3$}

There are a few computations in literature on the BV structures of
the Hochschild cohomology rings of some commutative algebras, see
for example, \cite{Yang2013}. As far as we know, there is no concrete
computation in non-commutative case. In this section, we use our
method to compute the BV structure of the Hochschild cohomology
rings of the group algebra $\mathbb{F}_3S_3$. The
associative ring structure has been determined by Siegel and
Witherspoon \cite{SW1999} using their cup product formula. So we
only need to compute the $\bigtriangleup$ operator and the Lie
bracket.

Let
$G=S_3=\langle a,b~|~a^3=1=b^2, bab=a^{-1}\rangle$. Choose the
conjugacy class representatives as $1,a,b$. The corresponding
centralizers are $H_1=G, H_2=\langle a\rangle$ and $H_3=\langle
b\rangle$. So $HH^*(\mathbb{F}_3S_3)\simeq H^*(S_3)\oplus
H^*(\langle a\rangle)\oplus H^*(\langle b\rangle)$. The ring
structures of $H^*(S_3)$, of $H^*(\langle a\rangle)$, and of
$H^*(\langle b\rangle)$ are well-known (see for example,
\cite{Evens}). $H^*(S_3)=\mathbb{F}_3[u,v]/(u^2)$, where $u$ and $v$
have degrees of $3$ and $4$, respectively. $H^*(\langle
a\rangle)=\mathbb{F}_3[w_1,w_2]/({w_1}^2)$, where $w_1$ and $w_2$
have degrees of $1$ and $2$, respectively. $H^*(\langle
b\rangle)=\mathbb{F}_3$, since $\mathbb{F}_3\langle b\rangle$ is
semisimple. Identify the elements $u, v$ with their images in
$HH^*(\mathbb{F}_3S_3)$ and denote by $X_1,X_2$ the images of the
elements (resp.) $w_1, w_2$ under the additive decomposition. Then
Siegel and Witherspoon proved in \cite{SW1999} the following
presentation for the Hochschild cohomology ring
$HH^*(\mathbb{F}_3S_3)$: $HH^*(\mathbb{F}_3S_3)$ is generated as an
algebra by elements $u,v,C_1=1+a+a^2,C_2=b(1+a+a^2),X_1,X_2$ of
degrees (resp.) $3,4,0,0,1$ and $2$, subject to the relations
$$uX_1=0,\quad vX_1=uX_2,\quad uC_2=0=vC_2,$$
$$
C_iX_j=0=C_iC_j (i,j\in \{1,2\}),\quad X_1X_2=uC_1,\quad
X_2^2=vC_1$$ in addition to the graded commutative relations.

\medskip
Our formula in Theorem \ref{bigtriangleup-operator} for
$\bigtriangleup$ operator is based on the normalized bar resolution.
However, the real computations of the Hochschild cohomology or the
group cohomology are based on the minimal projective resolutions. So
we need to construct comparison maps between the minimal projective
resolution and normalized bar resolution (by the same technique
introduced in Section 2), and then we can transfer our formula in
Theorem \ref{bigtriangleup-operator} to the minimal Hochschild
cohomology level. By Theorem \ref{bigtriangleup-operator}, the
operator $\bigtriangleup: HH^n(\mathbb{F}_3S_3)\longrightarrow
HH^{n-1}(\mathbb{F}_3S_3)$ restricts to the operators
$\widehat{\bigtriangleup}_b: H^n(\langle b\rangle)\longrightarrow
H^{n-1}(\langle b\rangle)$, $\widehat{\bigtriangleup}_a: H^n(\langle
a\rangle)\longrightarrow H^{n-1}(\langle a\rangle)$, and
$\widehat{\bigtriangleup}_1: H^n(S_3)\longrightarrow H^{n-1}(S_3)$.
Since $\mathbb{F}_3\langle b\rangle$ is semisimple and $H^*(\langle
b\rangle)$ is concentrated in degree zero,
$\widehat{\bigtriangleup}_b$ is trivial.

To compute $\widehat{\bigtriangleup}_a$, we first recall the minimal
projective resolution $P_a^*$ of the trivial $\mathbb{F}_3\langle
a\rangle$-module $\mathbb{F}_3$:
$$\cdots \longrightarrow \mathbb{F}_3\langle
a\rangle\stackrel{a-1}{\longrightarrow} \mathbb{F}_3\langle
a\rangle\stackrel{1+a+a^2}{\longrightarrow} \mathbb{F}_3\langle
a\rangle\stackrel{a-1}{\longrightarrow} \mathbb{F}_3\langle
a\rangle\stackrel{\epsilon}{\longrightarrow}
\mathbb{F}_3\longrightarrow 0,$$where the differential $\epsilon$ is
given by
$\epsilon(\lambda_1+\lambda_2a+\lambda_3a^2)=\lambda_1+\lambda_2+\lambda_3$,
and the differential $a-1$ means multiplying by $a-1$, etc.. There
is a setwise self-homotopy over $P_a^*$ as follows:
$$t_{-1}: \mathbb{F}_3\to \mathbb{F}_3\langle
a\rangle, \quad 1\mapsto 1,$$
$$t_0: \mathbb{F}_3\langle
a\rangle\to \mathbb{F}_3\langle a\rangle, \quad 1\mapsto 0, a\mapsto
1, a^2\mapsto 1+a,$$
$$t_1: \mathbb{F}_3\langle
a\rangle\to \mathbb{F}_3\langle a\rangle, \quad 1\mapsto 0, a\mapsto
0, a^2\mapsto 1,$$
$$t_2: \mathbb{F}_3\langle
a\rangle\to \mathbb{F}_3\langle a\rangle, \quad 1\mapsto 0, a\mapsto
1, a^2\mapsto 1+a,$$
$$t_3: \mathbb{F}_3\langle
a\rangle\to \mathbb{F}_3\langle a\rangle, \quad 1\mapsto 0, a\mapsto
0, a^2\mapsto 1,$$
$$\cdots \cdots \cdots$$
We also have the normalized bar resolution
$Bar_*(\mathbb{F}_3\langle a\rangle)\otimes_{\mathbb{F}_3\langle
a\rangle}\mathbb{F}_3$ of the trivial $\mathbb{F}_3\langle
a\rangle$-module $\mathbb{F}_3$, which is identified as the
following complex
$$\cdots \longrightarrow \mathbb{F}_3\langle
a\rangle\otimes \overline{\mathbb{F}_3\langle a\rangle}^{\otimes
n}\stackrel{d_n}{\longrightarrow} \cdots \longrightarrow
\mathbb{F}_3\langle a\rangle\otimes\overline{\mathbb{F}_3\langle
a\rangle}\stackrel{d_1}{\longrightarrow} \mathbb{F}_3\langle
a\rangle\stackrel{d_0}{\longrightarrow} \mathbb{F}_3\longrightarrow
0,$$where the differential is given by
$$d_0(g_0)=1  \quad (\mbox{for } g_0\in \langle
a\rangle)$$
 and (for $g_0\in \langle
a\rangle, g_1,\cdots, g_{n}\in \overline{\langle a\rangle}$)
$$\quad  \quad \quad  d_n(g_0, g_1,\cdots, g_{n})=g_0 g_1\otimes g_2\otimes\cdots
\otimes g_{n}+ \quad  \quad  \quad  \quad  \quad  \quad \quad \quad
\quad  \quad  \quad  \quad  \quad  \quad  \quad  \quad \quad \quad
\quad  \quad \quad  \quad  \quad  \quad  \quad  \quad \quad \quad
\quad \quad  \quad $$ $$\quad  \quad  \quad  \quad \quad \quad \quad
\quad \quad  \quad \quad  \quad  \quad \quad
\sum_{i=1}^{n-1}(-1)^ig_0\otimes \cdots \otimes g_ig_{i+1}\otimes
\cdots \otimes g_{n}+(-1)^{n}g_0\otimes g_1\otimes \cdots \otimes
g_{n-1}.$$ There is a setwise self-homotopy over $Bar_*(\mathbb{F}_3\langle a\rangle)\otimes_{\mathbb{F}_3\langle
a\rangle}\mathbb{F}_3$ as follows:
$$s_n: \mathbb{F}_3\langle
a\rangle\otimes \overline{\mathbb{F}_3\langle a\rangle}^{\otimes
n}\longrightarrow \mathbb{F}_3\langle a\rangle\otimes
\overline{\mathbb{F}_3\langle a\rangle}^{\otimes n+1},$$
$$g_0\otimes g_1\otimes \cdots \otimes g_{n}\longmapsto
1\otimes g_0\otimes g_1\otimes \cdots \otimes g_{n},$$ where $g_0\in
\langle a\rangle, g_1,\cdots, g_{n}\in \overline{\langle a\rangle}$.
Using $\{s_n\}$ and $\{t_n\}$ we get comparison maps $\Phi:
P_a^*\longrightarrow Bar_*(\mathbb{F}_3\langle
a\rangle)\otimes_{\mathbb{F}_3\langle a\rangle}\mathbb{F}_3$ and
$\Psi: Bar_*(\mathbb{F}_3\langle
a\rangle)\otimes_{\mathbb{F}_3\langle
a\rangle}\mathbb{F}_3\longrightarrow P_a^*$. We write down the maps
up to degree $4$ explicitly:
$$\Phi:
P_a^*\longrightarrow Bar_*(\mathbb{F}_3\langle
a\rangle)\otimes_{\mathbb{F}_3\langle a\rangle}\mathbb{F}_3$$
$$\Phi_{-1}=id:
\mathbb{F}_3\longrightarrow \mathbb{F}_3,$$
$$\Phi_{0}=id:
\mathbb{F}_3\langle a\rangle\longrightarrow \mathbb{F}_3\langle
a\rangle,$$
$$\Phi_{1}:
\mathbb{F}_3\langle a\rangle\longrightarrow \mathbb{F}_3\langle
a\rangle\otimes\overline{\mathbb{F}_3\langle a\rangle}, \quad
g\mapsto g\otimes a, \mbox{ for } g=1,a,a^2,$$
$$\Phi_{2}:
\mathbb{F}_3\langle a\rangle\longrightarrow \mathbb{F}_3\langle
a\rangle\otimes\overline{\mathbb{F}_3\langle a\rangle}^{\otimes 2},
\quad g\mapsto g\otimes a\otimes a + g\otimes a^2\otimes a, \mbox{
for } g=1,a,a^2,$$
$$\Phi_{3}:
\mathbb{F}_3\langle a\rangle\longrightarrow \mathbb{F}_3\langle
a\rangle\otimes\overline{\mathbb{F}_3\langle a\rangle}^{\otimes 3},
\quad g\mapsto g\otimes a\otimes a\otimes a + g\otimes a\otimes
a^2\otimes a, \mbox{ for } g=1,a,a^2,$$
$$\Phi_{4}:
\mathbb{F}_3\langle a\rangle\longrightarrow \mathbb{F}_3\langle
a\rangle\otimes\overline{\mathbb{F}_3\langle a\rangle}^{\otimes 4},
\quad g\mapsto g\otimes a\otimes a\otimes a\otimes a + g\otimes a\otimes a\otimes a^2\otimes a + $$
$$\quad  \quad  \quad  \quad  \quad  \quad \quad g\otimes a^2\otimes a\otimes
a\otimes a + g\otimes a^2\otimes a\otimes
a^2\otimes a, \mbox{ for } g=1,a,a^2;$$
$$\Psi:
Bar_*(\mathbb{F}_3\langle a\rangle)\otimes_{\mathbb{F}_3\langle
a\rangle}\mathbb{F}_3\longrightarrow P_a^*$$
$$\Psi_{-1}=id:
\mathbb{F}_3\longrightarrow \mathbb{F}_3,$$
$$\Psi_{0}=id:
\mathbb{F}_3\langle a\rangle\longrightarrow \mathbb{F}_3\langle
a\rangle,$$
$$\Psi_{1}:
\mathbb{F}_3\langle a\rangle\otimes\overline{\mathbb{F}_3\langle
a\rangle}\longrightarrow \mathbb{F}_3\langle a\rangle, \quad
g\otimes a\mapsto g, \quad g\otimes a^2\mapsto g(1+a),\mbox{ for }
g=1,a,a^2,$$
$$\Psi_{2}:
\mathbb{F}_3\langle a\rangle\otimes\overline{\mathbb{F}_3\langle
a\rangle}^{\otimes 2}\longrightarrow \mathbb{F}_3\langle a\rangle,
\quad g\otimes a\otimes a\mapsto 0, $$
$$g\otimes a\otimes a^2\mapsto g, \quad g\otimes a^2\otimes a\mapsto g, \quad g\otimes a^2\otimes a^2\mapsto g,
\mbox{ for } g=1,a,a^2,$$
$$\Psi_{3}:
\mathbb{F}_3\langle a\rangle\otimes\overline{\mathbb{F}_3\langle
a\rangle}^{\otimes 3}\longrightarrow \mathbb{F}_3\langle a\rangle,
\quad g\otimes a\otimes a\otimes a\mapsto 0, \quad g\otimes a\otimes
a\otimes a^2\mapsto 0,$$
$$g\otimes a\otimes
a^2\otimes a\mapsto 0, \quad g\otimes a\otimes a^2\otimes a^2\mapsto
0, \quad g\otimes a^2\otimes a\otimes a\mapsto 0, \quad $$
$$g\otimes a^2\otimes
a\otimes a^2\mapsto g(1+a), \quad g\otimes a^2\otimes a^2\otimes a\mapsto
g(1+a), \quad g\otimes a^2\otimes a^2\otimes a^2\mapsto ga, \mbox{ for }
g=1,a,a^2,$$
$$\Psi_{4}:
\mathbb{F}_3\langle a\rangle\otimes\overline{\mathbb{F}_3\langle
a\rangle}^{\otimes 4}\longrightarrow \mathbb{F}_3\langle a\rangle,
\quad g\otimes a\otimes a\otimes a\otimes a\mapsto 0, \quad g\otimes a\otimes a\otimes
a\otimes a^2\mapsto 0,$$
$$g\otimes a\otimes a\otimes
a^2\otimes a\mapsto 0, \quad g\otimes a\otimes a\otimes a^2\otimes a^2\mapsto
0, \quad g\otimes a\otimes a^2\otimes a\otimes a\mapsto 0, \quad $$
$$g\otimes a\otimes a^2\otimes
a\otimes a^2\mapsto g, \quad g\otimes a\otimes a^2\otimes a^2\otimes
a\mapsto g, \quad g\otimes a\otimes a^2\otimes a^2\otimes a^2\mapsto
g,$$
$$g\otimes a^2\otimes
a\otimes a\otimes a\mapsto 0, \quad g\otimes a^2\otimes a\otimes
a\otimes a^2\mapsto 0, \quad g\otimes a^2\otimes a\otimes a^2\otimes
a\mapsto 0,$$
$$g\otimes a^2\otimes a\otimes
a^2\otimes a^2\mapsto 0, \quad g\otimes a^2\otimes a^2\otimes
a\otimes a\mapsto 0, \quad g\otimes a^2\otimes a^2\otimes a\otimes
a^2\mapsto g,$$
$$g\otimes a^2\otimes a^2\otimes a^2\otimes a\mapsto
g, \quad g\otimes a^2\otimes a^2\otimes a^2\otimes a^2\mapsto 0,
\mbox{ for } g=1,a,a^2.$$

\medskip
We have the following commutative diagram:

\begin{center} {\setlength{\unitlength}{0.5cm}
\begin{picture}(27,7)
\put(13.2,3.5){$\mathbb{F}_3\langle
a\rangle$}\put(16.3,3.7){\vector(1,0){2}}\put(16.0,5.7){$\mathbb{F}_3$}\put(14.0,5.0){$w_1$}\put(14.3,4.2){\vector(1,1){1.5}}
\put(19,3.5){$\mathbb{F}_3\langle
a\rangle$}
\put(16.5,3.9){$a-1$}\put(22,3.9){$\epsilon$}\put(21,3.7){\vector(1,0){2}}\put(23.7,3.5){$\mathbb{F}_3$}
\put(25,3.7){\vector(1,0){2}}\put(27.5,3.5){$0$}\put(25,0.7){\vector(1,0){2}}\put(27.5,0.5){$0$}
\put(12.0,0.5){$\mathbb{F}_3\langle
a\rangle\otimes \overline{\mathbb{F}_3\langle a\rangle}$}\put(17.1,1.0){$d_1$}\put(16.3,0.7){\vector(1,0){2}}
\put(19,0.5){$\mathbb{F}_3\langle
a\rangle$}\put(22,1.0){$d_0$}\put(21,0.7){\vector(1,0){2}}\put(23.7,0.5){$\mathbb{F}_3$}
\put(13.7,1.5){\vector(0,1){1.5}}\put(12.8,2.0){$\Psi_{1}$}\put(14.4,2.0){$\Phi_{1}$}
\put(20.0,3.0){\vector(0,-1){1.5}}\put(18.6,2.0){$\Psi_{0}$}\put(14.3,3.0){\vector(0,-1){1.5}}\put(20.1,2.0){$\Phi_{0}$}
\put(19.5,1.5){\vector(0,1){1.5}}\put(23.8,2.0){$\|$}\put(23.8,1.5){$\|$}\put(23.8,2.5){$\|$}
\put(6.8,3.0){\vector(0,-1){1.5}}\put(6.3,1.5){\vector(0,1){1.5}}\put(5.3,2.0){$\Psi_{2}$}\put(6.9,2.0){$\Phi_{2}$}
\put(9.3,3.9){$1+a+a^2$}\put(10.3,1.0){$d_2$}\put(9.6,0.7){\vector(1,0){2}}\put(9.6,3.7){\vector(1,0){2}}
\put(6.0,3.5){$\mathbb{F}_3\langle
a\rangle$}\put(9.0,5.7){$\mathbb{F}_3$}\put(7.0,5.0){$w_2$}\put(7.3,4.2){\vector(1,1){1.5}}\put(4.5,0.5){$\mathbb{F}_3\langle
a\rangle\otimes \overline{\mathbb{F}_3\langle a\rangle}^{\otimes 2}$}
\put(2.4,3.9){$a-1$}\put(3.0,1.0){$d_3$}\put(2.2,0.7){\vector(1,0){2}}\put(2.2,3.7){\vector(1,0){2}}\put(0.0,3.9){$\ldots$}\put(0.0,1.0){$\ldots$}
\end{picture}
}\end{center} Clearly both the representatives of $w_1$ and $w_2$ in
the group cohomology $H^*(\langle
a\rangle)=\mathbb{F}_3[w_1,w_2]/({w_1}^2)$ can be chosen as
$\epsilon: \mathbb{F}_3\langle a\rangle \longrightarrow
\mathbb{F}_3, \lambda_1+\lambda_2a+\lambda_3a^2\mapsto
\lambda_1+\lambda_2+\lambda_3$. By abuse of notation, we have
$$\widehat{\bigtriangleup}_a(w_1)=\widehat{\bigtriangleup}_a(w_1\Psi_1)\circ \Phi_0,$$
$$\widehat{\bigtriangleup}_a(w_2)=\widehat{\bigtriangleup}_a(w_2\Psi_2)\circ \Phi_1.$$
A straightforward calculation shows that
$\widehat{\bigtriangleup}_a(w_1)=-1$ and
$\widehat{\bigtriangleup}_a(w_2)=0$. Similarly, we can get that
$$\widehat{\bigtriangleup}_a(w_1w_2)=
\widehat{\bigtriangleup}_a((w_1\Psi_1)(w_2\Psi_2))\circ\Phi_2,$$
$$\widehat{\bigtriangleup}_a(w_2^2)=\widehat{\bigtriangleup}_a((w_2\Psi_2)^2)\circ \Phi_3.$$
By direct computation, we have $(w_1\Psi_1)(w_2\Psi_2)\Phi_3(1)=1$ and $(w_2\Psi_2)^2)\Phi_4(1)=1$, which imply that both the representatives of $w_1w_2$ and $w_2^2$ in the group
cohomology $H^*(\langle a\rangle)=\mathbb{F}_3[w_1,w_2]/({w_1}^2)$
are given by $\epsilon$. So again a straightforward
calculation shows that $\widehat{\bigtriangleup}_a(w_1w_2)=-w_2$ and
$\widehat{\bigtriangleup}_a(w_2^2)=0$.

Next we compute $\widehat{\bigtriangleup}_1$. First of all, we need
to construct a minimal projective resolution $P_1^*$ of the trivial
$\mathbb{F}_3S_3$-module $\mathbb{F}_3$. Recall that the group
algebra $\mathbb{F}_3S_3$ can be identified as the
$\mathbb{F}_3$-algebra $A$ given by the following quiver and
relations:
$$\unitlength=1.00mm
\special{em:linewidth 0.4pt} \linethickness{0.4pt}
\begin{picture}(20.00,8.00)
\put(-3,3){$\stackrel{1}{\circ}$}
\put(12,3){$\stackrel{2}{\circ}$}\put(15,3){$,$}
\put(1,5){\vector(1,0){8.0}} \put(9,3){\vector(-1,0){8.0}}
\put(4,6){$\alpha$} \put(4,-1){$\beta$}
\put(20,3){$\alpha\beta\alpha=\beta\alpha\beta=0.$}
\end{picture}$$
Let $_AA=Ae_1\oplus Ae_2=\langle e_1, \alpha, \beta\alpha\rangle
\oplus \langle e_2, \beta, \alpha\beta\rangle$ be the decomposition
of the regular module into the indecomposable projective modules.
Then we have the following (Remind that all the computations take
place over $\mathbb{F}_3$):

$$e_1=-(1+b), e_2=-(1-b),$$
$$
\alpha=-a(1-a)(1+b)=-(1-b)a(1-a)=-(a-a^2+ab-a^2b),$$
$$\beta=-a(1-a)(1-b)=-(1+b)a(1-a)=-(a-a^2-ab+a^2b),$$
$$\beta\alpha=(1+b)(1+a+a^2)(1+b)=-(1+a+b+a^2+ab+a^2b),$$
$$\alpha\beta=(1-b)(1+a+a^2)(1-b)=-(1+a-b+a^2-ab-a^2b),
$$
$$
1=e_1+e_2,
a=e_1+e_2-\alpha-\beta-\alpha\beta-\beta\alpha,
b=e_1-e_2,$$
$$
a^2=e_1+e_2+\alpha+\beta-\alpha\beta-\beta\alpha,
ab=e_1-e_2-\alpha+\beta+\alpha\beta-\beta\alpha,$$
$$
a^2b=e_1-e_2+\alpha-\beta+\alpha\beta-\beta\alpha.$$
$Ae_1$ is the projective cover of the trivial module $\mathbb{F}_3$
since $a$ and $b$ act trivially on $Ae_1/rad(Ae_1)$; $Ae_2$ is the
projective cover of the sign module $sgn$ since $a$ (resp., $b$)
acts trivially (resp., by multiplying $-1$) on $Ae_2/rad(Ae_2)$. Now
it is easy to write down the minimal projective resolution $P_1^*$
of the trivial $\mathbb{F}_3S_3$-module $\mathbb{F}_3$ is as
follows:
$$\cdots \longrightarrow Ae_1\stackrel{\partial_4}{\longrightarrow}Ae_1\stackrel{\partial_3}{\longrightarrow}
Ae_2\stackrel{\partial_2}{\longrightarrow}
Ae_2\stackrel{\partial_1}{\longrightarrow}
Ae_1\stackrel{\partial_0}{\longrightarrow}
\mathbb{F}_3\longrightarrow 0,$$where the differential is given as
follows:
$$\partial_0: e_1\mapsto 1, \alpha\mapsto 0, \beta\alpha\mapsto 0,$$
$$\partial_1: e_2\mapsto \alpha, \beta\mapsto \beta\alpha, \alpha\beta\mapsto 0,$$
$$\partial_2: e_2\mapsto \alpha\beta, \beta\mapsto 0, \alpha\beta\mapsto 0,$$
$$\partial_3: e_1\mapsto \beta, \alpha\mapsto \alpha\beta, \beta\alpha\mapsto 0,$$
$$\partial_4: e_1\mapsto \beta\alpha, \alpha\mapsto 0, \beta\alpha\mapsto 0,$$
$$\partial_5: e_2\mapsto \alpha, \beta\mapsto \beta\alpha, \alpha\beta\mapsto 0,$$
$$\partial_6: e_2\mapsto \alpha\beta, \beta\mapsto 0, \alpha\beta\mapsto 0,$$
$$\partial_7: e_1\mapsto \beta, \alpha\mapsto \alpha\beta, \beta\alpha\mapsto 0,$$
$$\partial_8: e_1\mapsto \beta\alpha, \alpha\mapsto 0, \beta\alpha\mapsto 0,$$
$$\cdots \cdots \cdots$$
Using the Lowey diagram structures of $Ae_1$ and $Ae_2$, we can
easily construct a setwise self-homotopy over $P_1^*$ as follows:
$$t_{-1}: \mathbb{F}_3\to Ae_1, \quad 1\mapsto e_1,$$
$$t_0: Ae_1\to Ae_2, \quad e_1\mapsto 0, \alpha\mapsto
e_2, \beta\alpha\mapsto \beta,$$
$$t_1: Ae_2\to Ae_2, \quad e_2\mapsto 0, \beta\mapsto
0, \alpha\beta\mapsto e_2,$$
$$t_2: Ae_2\to Ae_1, \quad e_2\mapsto 0, \beta\mapsto
e_1, \alpha\beta\mapsto \alpha,$$
$$t_3: Ae_1\to Ae_1, \quad e_1\mapsto 0, \alpha\mapsto
0, \beta\alpha\mapsto e_1,$$
$$t_4: Ae_1\to Ae_2, \quad e_1\mapsto 0, \alpha\mapsto
e_2, \beta\alpha\mapsto \beta,$$
$$t_5: Ae_2\to Ae_2, \quad e_2\mapsto 0, \beta\mapsto
0, \alpha\beta\mapsto e_2,$$
$$t_6: Ae_2\to Ae_1, \quad e_2\mapsto 0, \beta\mapsto
e_1, \alpha\beta\mapsto \alpha,$$
$$t_7: Ae_1\to Ae_1, \quad e_1\mapsto 0, \alpha\mapsto
0, \beta\alpha\mapsto e_1,$$
$$t_8: Ae_1\to Ae_2, \quad e_1\mapsto 0, \alpha\mapsto
e_2, \beta\alpha\mapsto \beta,$$
$$\cdots \cdots \cdots$$
We also have the normalized bar resolution
$Bar_*(A)\otimes_{A}\mathbb{F}_3$ of the trivial $A$-module $\mathbb{F}_3$, which is identified as the
following complex
$$\cdots \longrightarrow A\otimes \overline{A}^{\otimes
n}\stackrel{d_n}{\longrightarrow} \cdots \longrightarrow
A\otimes\overline{A}\stackrel{d_1}{\longrightarrow} A\stackrel{d_0}{\longrightarrow} \mathbb{F}_3\longrightarrow
0,$$where the differential is given by
$$d_0(g_0)=1  \quad (\mbox{for } g_0\in S_3)$$
 and (for $g_0\in S_3, g_1,\cdots, g_{n}\in \overline{S_3}$)
$$\quad  \quad \quad  d_n(g_0, g_1,\cdots, g_{n})=g_0 g_1\otimes g_2\otimes\cdots
\otimes g_{n}+ \quad  \quad  \quad  \quad  \quad  \quad \quad \quad
\quad  \quad  \quad  \quad  \quad  \quad  \quad  \quad \quad \quad
\quad  \quad \quad  \quad  \quad  \quad  \quad  \quad \quad \quad
\quad \quad  \quad $$ $$\quad  \quad  \quad  \quad \quad \quad \quad
\quad \quad  \quad \quad  \quad  \quad \quad
\sum_{i=1}^{n-1}(-1)^ig_0\otimes \cdots \otimes g_ig_{i+1}\otimes
\cdots \otimes g_{n}+(-1)^{n}g_0\otimes g_1\otimes \cdots \otimes
g_{n-1}.$$ There is a setwise self-homotopy over $Bar_*(A)\otimes_{A}\mathbb{F}_3$ as follows:
$$s_n: A\otimes \overline{A}^{\otimes
n}\longrightarrow A\otimes
\overline{A}^{\otimes n+1},$$
$$g_0\otimes g_1\otimes \cdots \otimes g_{n}\longmapsto
1\otimes g_0\otimes g_1\otimes \cdots \otimes g_{n},$$ where $g_0\in
S_3, g_1,\cdots, g_{n}\in \overline{S_3}$.
As before, we want to use $\{s_n\}$ and $\{t_n\}$ to get comparison maps $\Lambda:
P_1^*\longrightarrow Bar_*(A)\otimes_{A}\mathbb{F}_3$ and
$\Theta: Bar_*(A)\otimes_{A}\mathbb{F}_3\longrightarrow P_1^*$. Here the situation is a bit different, since $P_1^*$
is not a free resolution. However, if we replace $s_n(x)$ by $\widetilde{s_n}(x)=e_1s_n(e_1x)+e_2s_n(e_2x)$, then the method introduced in
Section 2 still works. We write down the comparison maps
up to degree $8$ explicitly:
$$\Lambda:
P_1^*\longrightarrow Bar_*(A)\otimes_{A}\mathbb{F}_3$$
$$\Lambda_{-1}=id:
\mathbb{F}_3\longrightarrow \mathbb{F}_3,$$
$$\Lambda_{0}:
Ae_1\longrightarrow A,$$
$$ae_1\mapsto ae_1,$$
$$\Lambda_{1}:
Ae_2\longrightarrow A\otimes\overline{A},$$
$$
e_2\mapsto e_2\otimes \alpha,
\beta\mapsto \beta\otimes \alpha, \alpha\beta\mapsto \alpha\beta\otimes \alpha,$$
$$\Lambda_{2}:
Ae_2\longrightarrow A\otimes\overline{A}^{\otimes 2},$$
$$
e_2\mapsto e_2\otimes \alpha\beta\otimes \alpha,
\beta\mapsto \beta\otimes \alpha\beta\otimes \alpha,
\alpha\beta\mapsto \alpha\beta\otimes \alpha\beta\otimes \alpha,$$
$$\Lambda_{3}:
Ae_1\longrightarrow A\otimes\overline{A}^{\otimes 3},$$
$$
e_1\mapsto e_1\otimes\beta\otimes \alpha\beta\otimes \alpha,
\alpha\mapsto \alpha\otimes\beta\otimes \alpha\beta\otimes \alpha,
\beta\alpha\mapsto \beta\alpha\otimes\beta\otimes \alpha\beta\otimes \alpha,$$
$$\Lambda_{4}:
Ae_1\longrightarrow A\otimes\overline{A}^{\otimes 4},$$
$$
e_1\mapsto e_1\otimes\beta\alpha\otimes\beta\otimes
\alpha\beta\otimes \alpha, \alpha\mapsto
\alpha\otimes\beta\alpha\otimes\beta\otimes \alpha\beta\otimes
\alpha, \beta\alpha\mapsto
\beta\alpha\otimes\beta\alpha\otimes\beta\otimes \alpha\beta\otimes
\alpha,$$
$$\Lambda_{5}:
Ae_2\longrightarrow A\otimes\overline{A}^{\otimes 5},$$
$$
e_2\mapsto e_2\otimes \alpha\otimes\beta\alpha\otimes\beta\otimes
\alpha\beta\otimes \alpha, \beta\mapsto \beta\otimes
\alpha\otimes\beta\alpha\otimes\beta\otimes \alpha\beta\otimes
\alpha,$$
$$\alpha\beta\mapsto \alpha\beta\otimes
\alpha\otimes\beta\alpha\otimes\beta\otimes \alpha\beta\otimes
\alpha,$$
$$\Lambda_{6}:
Ae_2\longrightarrow A\otimes\overline{A}^{\otimes 6},$$
$$
e_2\mapsto e_2\otimes \alpha\beta\otimes
\alpha\otimes\beta\alpha\otimes\beta\otimes \alpha\beta\otimes
\alpha, \beta\mapsto \beta\otimes \alpha\beta\otimes
\alpha\otimes\beta\alpha\otimes\beta\otimes \alpha\beta\otimes
\alpha,$$
$$\alpha\beta\mapsto \alpha\beta\otimes \alpha\beta\otimes
\alpha\otimes\beta\alpha\otimes\beta\otimes \alpha\beta\otimes
\alpha,$$
$$\Lambda_{7}:
Ae_1\longrightarrow A\otimes\overline{A}^{\otimes 7},$$
$$
e_1\mapsto e_1\otimes\beta\otimes \alpha\beta\otimes
\alpha\otimes\beta\alpha\otimes\beta\otimes \alpha\beta\otimes
\alpha, \alpha\mapsto \alpha\otimes\beta\otimes \alpha\beta\otimes
\alpha\otimes\beta\alpha\otimes\beta\otimes \alpha\beta\otimes
\alpha,$$
$$\beta\alpha\mapsto \beta\alpha\otimes\beta\otimes \alpha\beta\otimes
\alpha\otimes\beta\alpha\otimes\beta\otimes \alpha\beta\otimes
\alpha,$$
$$\Lambda_{8}:
Ae_1\longrightarrow A\otimes\overline{A}^{\otimes 8},$$
$$
e_1\mapsto e_1\otimes\beta\alpha\otimes\beta\otimes
\alpha\beta\otimes \alpha\otimes\beta\alpha\otimes\beta\otimes
\alpha\beta\otimes \alpha, \alpha\mapsto
\alpha\otimes\beta\alpha\otimes\beta\otimes \alpha\beta\otimes
\alpha\otimes\beta\alpha\otimes\beta\otimes \alpha\beta\otimes
\alpha,$$
$$\beta\alpha\mapsto
\beta\alpha\otimes\beta\alpha\otimes\beta\otimes \alpha\beta\otimes
\alpha\otimes\beta\alpha\otimes\beta\otimes \alpha\beta\otimes
\alpha;$$
$$\Theta:
Bar_*(A)\otimes_{A}\mathbb{F}_3\longrightarrow P_1^*$$
$$\Theta_{-1}=id:
\mathbb{F}_3\longrightarrow \mathbb{F}_3,$$
$$\Theta_{0}:
A\longrightarrow Ae_1,$$
$$a\mapsto ae_1,\mbox{ for } a\in A,$$
$$\Theta_{1}:
A\otimes\overline{A}\longrightarrow Ae_2,$$
$$ 1\otimes g_1\mapsto
-e_2-\beta, 0, e_2-\beta, -e_2-\beta, e_2-\beta,\mbox{ for }
g_1=a,b,a^2,ab,a^2b,$$
$$\Theta_{2}:
A\otimes\overline{A}^{\otimes 2}\longrightarrow Ae_2,$$
$$1\otimes g_1\otimes g_2\mapsto -e_2, 0, 0, 0, e_2 \mbox{ for } g_1=a,b,a^2,ab,a^2b,$$
$$ \mbox{ and where }
\Theta_{1}(1\otimes g_2)=-e_2-\beta,$$
$$1\otimes g_1\otimes g_2\mapsto 0, 0, e_2, -e_2, 0 \mbox{ for } g_1=a,b,a^2,ab,a^2b,$$
$$ \mbox{ and where }
\Theta_{1}(1\otimes g_2)=e_2-\beta,$$
$$1\otimes g_1\otimes g_2\mapsto 0, \mbox{ for any other case,} $$
$$\Theta_{3}:
A\otimes\overline{A}^{\otimes 3}\longrightarrow Ae_1,$$
$$1\otimes
g_1\otimes g_2\otimes g_3\mapsto -e_1-\alpha, 0, e_1-\alpha,
e_1+\alpha, -e_1+\alpha,  \mbox{ for } g_1=a,b,a^2,ab,a^2b,$$
$$\mbox{ and where } \Theta_{2}(1\otimes g_2\otimes g_3)=e_2,$$
$$1\otimes
g_1\otimes g_2\otimes g_3\mapsto e_1+\alpha, 0, -e_1+\alpha,
-e_1-\alpha, e_1-\alpha,  \mbox{ for } g_1=a,b,a^2,ab,a^2b,$$
$$\mbox{ and where } \Theta_{2}(1\otimes g_2\otimes g_3)=-e_2,$$
$$1\otimes g_1\otimes g_2\otimes g_3\mapsto 0, \mbox{ for any other case,} $$
$$\Theta_{4}: A\otimes\overline{A}^{\otimes
4}\longrightarrow Ae_1,$$
$$1\otimes g_1\otimes g_2\otimes g_3\otimes
g_4\mapsto e_1, 0, 0, 0, e_1 \mbox{ for } g_1=a,b,a^2,ab,a^2b,$$
$$ \mbox{ and where }
\Theta_{3}(1\otimes g_2\otimes g_3\otimes g_4)=e_1+\alpha,$$
$$1\otimes g_1\otimes g_2\otimes g_3\otimes
g_4\mapsto 0, 0, e_1, e_1, 0 \mbox{ for } g_1=a,b,a^2,ab,a^2b,$$
$$ \mbox{ and where }
\Theta_{3}(1\otimes g_2\otimes g_3\otimes g_4)=e_1-\alpha,$$
$$1\otimes g_1\otimes g_2\otimes g_3\otimes
g_4\mapsto -e_1, 0, 0, 0, -e_1 \mbox{ for } g_1=a,b,a^2,ab,a^2b,$$
$$ \mbox{ and where }
\Theta_{3}(1\otimes g_2\otimes g_3\otimes g_4)=-e_1-\alpha,$$
$$1\otimes g_1\otimes g_2\otimes g_3\otimes
g_4\mapsto 0, 0, -e_1, -e_1, 0 \mbox{ for } g_1=a,b,a^2,ab,a^2b,$$
$$ \mbox{ and where }
\Theta_{3}(1\otimes g_2\otimes g_3\otimes g_4)=-e_1+\alpha,$$
$$1\otimes g_1\otimes g_2\otimes g_3\otimes g_4\mapsto 0, \mbox{ for any other case,} $$
$$\Theta_{5}: A\otimes\overline{A}^{\otimes
5}\longrightarrow Ae_2,$$
$$1\otimes g_1\otimes g_2\otimes g_3\otimes
g_4\otimes g_5\mapsto -e_2-\beta, 0, e_2-\beta, -e_2-\beta,
e_2-\beta,  \mbox{ for } g_1=a,b,a^2,ab,a^2b,$$
$$ \mbox{ and where }
\Theta_{4}(1\otimes g_2\otimes g_3\otimes g_4\otimes g_5)=e_1,$$
$$1\otimes g_1\otimes g_2\otimes g_3\otimes
g_4\otimes g_5\mapsto e_2+\beta, 0, -e_2+\beta, e_2+\beta,
-e_2+\beta,  \mbox{ for } g_1=a,b,a^2,ab,a^2b,$$
$$ \mbox{ and where }
\Theta_{4}(1\otimes g_2\otimes g_3\otimes g_4\otimes g_5)=-e_1,$$
$$1\otimes g_1\otimes g_2\otimes g_3\otimes g_4\otimes g_5\mapsto 0, \mbox{ for any other case,} $$
$$\Theta_{6}: A\otimes\overline{A}^{\otimes
6}\longrightarrow Ae_2,$$
$$1\otimes g_1\otimes g_2\otimes g_3\otimes
g_4\otimes g_5\otimes g_6\mapsto e_2, 0, 0, 0, -e_2 \mbox{ for }
g_1=a,b,a^2,ab,a^2b,$$
$$ \mbox{ and where }
\Theta_{5}(1\otimes g_2\otimes g_3\otimes g_4\otimes g_5\otimes
g_6)=e_2+\beta,$$
$$1\otimes g_1\otimes g_2\otimes g_3\otimes
g_4\otimes g_5\otimes g_6\mapsto 0, 0, e_2, -e_2, 0 \mbox{ for }
g_1=a,b,a^2,ab,a^2b,$$
$$ \mbox{ and where }
\Theta_{5}(1\otimes g_2\otimes g_3\otimes g_4\otimes g_5\otimes
g_6)=e_2-\beta,$$
$$1\otimes g_1\otimes g_2\otimes g_3\otimes
g_4\otimes g_5\otimes g_6\mapsto -e_2, 0, 0, 0, e_2 \mbox{ for }
g_1=a,b,a^2,ab,a^2b,$$
$$ \mbox{ and where }
\Theta_{5}(1\otimes g_2\otimes g_3\otimes g_4\otimes g_5\otimes
g_6)=-e_2-\beta,$$
$$1\otimes g_1\otimes g_2\otimes g_3\otimes
g_4\otimes g_5\otimes g_6\mapsto 0, 0, -e_2, e_2, 0 \mbox{ for }
g_1=a,b,a^2,ab,a^2b,$$
$$ \mbox{ and where }
\Theta_{5}(1\otimes g_2\otimes g_3\otimes g_4\otimes g_5\otimes
g_6)=-e_2+\beta,$$
$$1\otimes g_1\otimes g_2\otimes g_3\otimes
g_4\otimes g_5\otimes g_6\mapsto 0, \mbox{ for any other case,} $$
$$\Theta_{7}: A\otimes\overline{A}^{\otimes
7}\longrightarrow Ae_1,$$
$$1\otimes g_1\otimes g_2\otimes g_3\otimes
g_4\otimes g_5\otimes g_6\otimes g_7\mapsto -e_1-\alpha, 0,
e_1-\alpha, e_1+\alpha, -e_1+\alpha \mbox{ for }
g_1=a,b,a^2,ab,a^2b,$$
$$ \mbox{ and where }
\Theta_{6}(1\otimes g_2\otimes g_3\otimes g_4\otimes g_5\otimes
g_6\otimes g_7)=e_2,$$
$$1\otimes g_1\otimes g_2\otimes g_3\otimes
g_4\otimes g_5\otimes g_6\otimes g_7\mapsto e_1+\alpha, 0,
-e_1+\alpha, -e_1-\alpha, e_1-\alpha \mbox{ for }
g_1=a,b,a^2,ab,a^2b,$$
$$ \mbox{ and where }
\Theta_{6}(1\otimes g_2\otimes g_3\otimes g_4\otimes g_5\otimes
g_6\otimes g_7)=-e_2,$$
$$1\otimes g_1\otimes g_2\otimes g_3\otimes
g_4\otimes g_5\otimes g_6\otimes g_7\mapsto 0, \mbox{ for any other
case,}
$$$$\Theta_{8}: A\otimes\overline{A}^{\otimes 8}\longrightarrow
Ae_1,$$
$$1\otimes g_1\otimes g_2\otimes g_3\otimes
g_4\otimes g_5\otimes g_6\otimes g_7\otimes g_8\mapsto e_1, 0, 0, 0,
e_1 \mbox{ for } g_1=a,b,a^2,ab,a^2b,$$
$$ \mbox{ and where }
\Theta_{7}(1\otimes g_2\otimes g_3\otimes g_4\otimes g_5\otimes
g_6\otimes g_7\otimes g_8)=e_1+\alpha,$$
$$1\otimes g_1\otimes g_2\otimes g_3\otimes
g_4\otimes g_5\otimes g_6\otimes g_7\otimes g_8\mapsto 0, 0, e_1,
e_1, 0 \mbox{ for } g_1=a,b,a^2,ab,a^2b,$$
$$ \mbox{ and where }
\Theta_{7}(1\otimes g_2\otimes g_3\otimes g_4\otimes g_5\otimes
g_6\otimes g_7\otimes g_8)=e_1-\alpha,$$
$$1\otimes g_1\otimes g_2\otimes g_3\otimes
g_4\otimes g_5\otimes g_6\otimes g_7\otimes g_8\mapsto -e_1, 0, 0, 0,
-e_1 \mbox{ for } g_1=a,b,a^2,ab,a^2b,$$
$$ \mbox{ and where }
\Theta_{7}(1\otimes g_2\otimes g_3\otimes g_4\otimes g_5\otimes
g_6\otimes g_7\otimes g_8)=-e_1-\alpha,$$
$$1\otimes g_1\otimes g_2\otimes g_3\otimes
g_4\otimes g_5\otimes g_6\otimes g_7\otimes g_8\mapsto 0, 0, -e_1,
-e_1, 0 \mbox{ for } g_1=a,b,a^2,ab,a^2b,$$
$$ \mbox{ and where }
\Theta_{7}(1\otimes g_2\otimes g_3\otimes g_4\otimes g_5\otimes
g_6\otimes g_7\otimes g_8)=-e_1+\alpha,$$
$$1\otimes g_1\otimes g_2\otimes g_3\otimes
g_4\otimes g_5\otimes g_6\otimes g_7\otimes g_8\mapsto 0, \mbox{ for
any other case.} $$ Note that both the representatives of $u$ and
$v$ in the group cohomology $H^*(S_3)=\mathbb{F}_3[u,v]/(u^2)$ can
be chosen as $Ae_1 \longrightarrow \mathbb{F}_3, e_1\mapsto 1,
\alpha\mapsto 0, \beta\alpha\mapsto 0$. Since $|u|=3$ and $|v|=4$,
we have $\widehat{\bigtriangleup}_1(u)=0$ and
$\widehat{\bigtriangleup}_1(v)=\lambda u$ for some $\lambda\in
\mathbb{F}_3$. We have
$$\widehat{\bigtriangleup}_1(v)=\widehat{\bigtriangleup}_1(v\Theta_4)\circ \Lambda_3,$$
and $\widehat{\bigtriangleup}_1(v\Theta_4)$ can be computed by our
formula in Theorem \ref{bigtriangleup-operator}. Since
$$\widehat{\bigtriangleup}_1(v\Theta_4)\circ
\Lambda_3(e_1)=e_1\widehat{\bigtriangleup}_1(v\Theta_4)(\beta\otimes
\alpha\beta\otimes
\alpha)$$
$$=-\widehat{\bigtriangleup}_1(v\Theta_4)(\beta\otimes
\alpha\beta\otimes
\alpha)-b\widehat{\bigtriangleup}_1(v\Theta_4)(\beta\otimes
\alpha\beta\otimes \alpha)$$
$$=-2\widehat{\bigtriangleup}_1(v\Theta_4)(\beta\otimes
\alpha\beta\otimes
\alpha)=\widehat{\bigtriangleup}_1(v\Theta_4)(\beta\otimes
\alpha\beta\otimes \alpha).$$ By a MAPLE calculation (see
\cite{Maple}: A MAPLE program for computing
$\widehat{\bigtriangleup}_1.$), we obtain that
$\widehat{\bigtriangleup}_1(v\Theta_4)(\beta\otimes
\alpha\beta\otimes \alpha)=0$, and therefore
$\widehat{\bigtriangleup}_1(v)=0$.
 Since $|uv|=7$ and $|v^2|=8$,
we have $\widehat{\bigtriangleup}_1(uv)=\mu u^2$ and
$\widehat{\bigtriangleup}_1(v^2)=\mu' uv$ for some $\mu, \mu'\in
\mathbb{F}_3$. Since $u^2=0, \widehat{\bigtriangleup}_1(uv)=0$, and
we only need to compute $\widehat{\bigtriangleup}_1(v^2)$. The
representative of $v^2$ in the group cohomology
$H^*(S_3)=\mathbb{F}_3[u,v]/(u^2)$ can also be chosen (up to a sign) as $Ae_1
\longrightarrow \mathbb{F}_3, e_1\mapsto 1, \alpha\mapsto 0,
\beta\alpha\mapsto 0$.  We have
$$\widehat{\bigtriangleup}_1(v^2)=\widehat{\bigtriangleup}_1(v^2\Theta_8)\circ \Lambda_7,$$
 $$\widehat{\bigtriangleup}_1(v^2\Theta_8)\circ
\Lambda_7(e_1)=e_1\widehat{\bigtriangleup}_1(v^2\Theta_7)(\beta\otimes
\alpha\beta\otimes \alpha\otimes\beta\alpha\otimes\beta\otimes
\alpha\beta\otimes \alpha)$$
$$=\widehat{\bigtriangleup}_1(v^2\Theta_7)(\beta\otimes
\alpha\beta\otimes \alpha\otimes\beta\alpha\otimes\beta\otimes
\alpha\beta\otimes \alpha).$$ Similarly by a MAPLE
calculation (see\cite{Maple}: A MAPLE program for computing
$\widehat{\bigtriangleup}_1$), we obtain that
$\widehat{\bigtriangleup}_1(v^2\Theta_8)\circ \Lambda_7(e_1)=0$, and
therefore $\widehat{\bigtriangleup}_1(v^2)=0$.

Finally, based on the above computations, we deal with the Lie
brackets. Since we have the following Possion rule: $[\alpha\cup
\beta, \gamma]=[\alpha, \gamma]\cup \beta +
(-1)^{|\alpha|(|\gamma|-1)}\alpha\cup[\beta, \gamma]$, it suffices
to write down the Lie brackets between generators in
$HH^*(\mathbb{F}_3S_3)$. Recall that $HH^*(\mathbb{F}_3S_3)$ is
generated as an algebra by elements
$u,v,C_1=1+a+a^2,C_2=b(1+a+a^2),X_1,X_2$ of degrees (resp.)
$3,4,0,0,1$ and $2$, subject to the relations
$$uX_1=0,\quad vX_1=uX_2,\quad uC_2=0=vC_2,$$
$$
C_iX_j=0=C_iC_j (i,j\in \{1,2\}),\quad X_1X_2=uC_1,\quad
X_2^2=vC_1$$ in addition to the graded commutative relations. Using
the formulas (Here we omit the sign $\cup$ in the equation)
$$[\alpha,\beta]=-(-1)^{(|\alpha|-1)|\beta|}(\bigtriangleup(\alpha\beta)-\bigtriangleup(\alpha)\beta-(-1)^{|\alpha|}\alpha\bigtriangleup(\beta))$$
 and
$$[\alpha,\beta]=-(-1)^{(|\alpha|-1)(|\beta|-1)}[\beta,\alpha],$$ we do the concrete
computations as follows:
$$[u,u]=0, [u,v]=\widehat{\bigtriangleup}_1(uv)-\widehat{\bigtriangleup}_1(u)v+u\widehat{\bigtriangleup}_1(v)=0, [v,u]=0,$$
$$[u,C_1]=-(\bigtriangleup(u C_1)-\bigtriangleup(u) C_1+u\bigtriangleup(C_1))=-\widehat{\bigtriangleup}_a(X_1X_2)=X_2,
[C_1,u]=-X_2,$$
$$[u,C_2]=-(\bigtriangleup(u C_2)-\bigtriangleup(u) C_2+u\bigtriangleup(C_2))=0, [C_2,u]=0,$$
$$[u,X_1]=-(\bigtriangleup(u X_1)-\bigtriangleup(u) X_1+u\bigtriangleup(X_1))=u,[X_1,u]=-u,$$
$$[u,X_2]=-(\bigtriangleup(u X_2)-\bigtriangleup(u) X_2+u\bigtriangleup(X_2))=-\widehat{\bigtriangleup}(uX_2)=0^c\footnotetext[3]{$uX_2$ is an element of degree $5$,
under the additive decompositon, it correponds to an element in
$H^*(\langle a\rangle)$ and has the form $\pm w_1{w_2}^2$. It
follows from the formula in the last paragragh of Section 8 that
$\widehat{\bigtriangleup}(uX_2)=\widehat{\bigtriangleup}_a(\pm
w_1{w_2}^2)=0.$},[X_2,u]=0,$$
$$[v,v]=-(\bigtriangleup(v^2)-\bigtriangleup(v)v-v\bigtriangleup(v))=0,$$
$$[v,C_1]=-(\bigtriangleup(vC_1)-\bigtriangleup(v) C_1+v\bigtriangleup(C_1))=-\widehat{\bigtriangleup}_a(X_2^2)=0, [C_1,v]=0,$$
$$[v,C_2]=-(\bigtriangleup(v C_2)-\bigtriangleup(v) C_2+v\bigtriangleup(C_2))=0,[C_2,v]=0,$$
$$[v,X_1]=\bigtriangleup(v X_1)-\bigtriangleup(v) X_1-v\bigtriangleup(X_1)=-v,[X_1,v]=v,$$
$$[v,X_2]=-(\bigtriangleup(v X_2)-\bigtriangleup(v) X_2-v\bigtriangleup(X_2))=\pm\widehat{\bigtriangleup}_a(X_2^3)=0,[X_2,v]=0,$$
$$[C_1,C_1]=[C_1,C_2]=[C_2,C_1]=[C_2,C_2]=0,$$
$$[C_1,X_1]=\bigtriangleup(C_1 X_1)-\bigtriangleup(C_1)X_1 -C_1\bigtriangleup(X_1)=C_1,[X_1,C_1]=-C_1,$$
$$[C_1,X_2]=-(\bigtriangleup(C_1 X_2)-\bigtriangleup(C_1) X_2-C_1\bigtriangleup(X_2))=0,[X_2,C_1]=0,$$
$$[C_2,X_1]=\bigtriangleup(C_2X_1)-\bigtriangleup(C_2) X_1-C_2\bigtriangleup(X_1)=C_2,[X_1,C_2]=-C_2,$$
$$[C_2,X_2]=-(\bigtriangleup(C_2X_2)-\bigtriangleup(C_2) X_2-C_2\bigtriangleup(X_2))=0,[X_2,C_2]=0,$$
$$[X_1,X_1]=0,[X_1,X_2]=-(\bigtriangleup(X_1X_2)-\bigtriangleup(X_1)X_2+X_1\bigtriangleup(X_2))=0,[X_1,X_2]=0,$$
$$[X_2,X_2]=-(\bigtriangleup(X_2^2)-\bigtriangleup(X_2)X_2-X_2\bigtriangleup(X_2))=0.$$

\begin{Rem} By a recent result of Menichi (see \cite[p. 321]{Menichi2011}), the Lie bracket of the group cohomology $H^*(G)$ for a finite group $G$ must be trivial. The above computation shows that this is indeed the case for $H^*(S_3)=\mathbb{F}_3[u,v]/(u^2)$. Note that to verify $[v,v]=0$, we have used the MAPLE program in \cite{Maple}.

\end{Rem}

\begin{Rem} Observe in the above example that the generators of $HH^*(\mathbb{F}_3S_3)$ are ``multiplicative closed" under Lie bracket: the Lie bracket $[\alpha,\beta]$ of two generators $\alpha$ and $\beta$ is a scalar multiple of another generator. Also if $[\alpha,\beta]\neq 0$, then $[\alpha,\beta]$ is equal to $-[\beta,\alpha]$.

\end{Rem}

\end{document}